\documentclass[12pt]{article}

\usepackage[english]{babel} 		%choix de la langue
\usepackage{amsmath}
\usepackage{array}
\usepackage{amssymb}
\usepackage{enumerate,epsfig}		%pour les figures
\usepackage{graphicx}   		%standard LaTeX graphics tool 
\newtheorem{Th}{Theorem}[section]
\newtheorem{Prop}[Th]{Proposition}
\newtheorem{Lemma}[Th]{Lemma}

\newtheorem{Cor}[Th]{Corollary}
\newtheorem{Rk}[Th]{Remark}
\makeatletter
\@addtoreset{equation}{section}
\makeatother

\newcommand{\cqfd}{{\unskip\kern 6pt\penalty 500
\raise -2pt\hbox{\vrule\vbox to 6pt{\hrule width 6pt
\vfill\hrule}\vrule}\par}}
\def\del{\partial}			%le d rond
\def\eps{\varepsilon} 			%le petit epsilon de l'analyse
\def\N{{\mathbb N}} 			%ensemble des entiers naturels
\def\r{{\mathbb R}} 			%ensemble des r\'eels
\newcommand\R{{\mathbb R}}
\def\1{{1\hspace{-0.9mm}{\rm l}}}	%fonction indicatrice
\def\ds{\displaystyle}			%pour l'affichage en grand des symboles
\def\tw{\mathop{{\rightharpoonup}}\limits} %convergence faible
\newcommand{\ft}[1]		
{\left(\begin{array}{c} #1 \end{array} \right)}
\newcommand{\coord}[2]		%vecteur colonne de r^2
{\ft{#1\\#2}}
\newcommand{\coords}[3]		%vecteur colonne de r^3
{\ft{#1\\#2\\#3}}
\newcommand{\ligne}[3]	%vecteur ligne de r^n
{\left(\begin{array}{ccccc} #1 & #2 & \dots & #3 \end{array} \right)}
\newcommand{\colonne}[3]	%vecteur colonne de r^n
{\left(\begin{array}{c} #1\\#2\\ \vdots\\#3 \end{array} \right)}
\newcommand{\colonnehaute}[3]	%vecteur colonne de r^n
{\left(\begin{array}{c} #1\\#2\\ \vdots\\\vdots\\#3 \end{array} \right)}
\newcommand{\deter}[4]		%determinant 2x2
{\left|\begin{array}{cc} #1 & #2\\#3 & #4 \end{array}
\right|}
\newcommand{\mat}[4]		%matrice 2x2
{\left(\begin{array}{cc} #1 & #2\\#3 & #4 \end{array}
\right)}
\newcommand{\determ}[9]		%determinant 3x3
{\left|\begin{array}{ccc} #1 & #2 & #3\\#4 & #5
& #6\\ #7 & #8 & #9\end{array}
\right|}
\newcommand{\matr}[9]		%matrice 3x3
{\left(\begin{array}{ccc} #1 & #2 & #3\\#4 & #5
& #6\\ #7 & #8 & #9 \end{array}
\right)}
\newcommand{\determinant}[9]	%determinant nxn
{\left|\begin{array}{cccc} #1 & #2 &\dots &  #3\\#4
& #5
& \dots & #6\\ #7 & #8 & \dots & #9\end{array}
\right|}
\newcommand{\matrice}[9]	%matrice nxn
{\left(\begin{array}{cccc} #1 & #2 &\dots &  #3\\#4
& #5
& \dots & #6\\ \vdots & \vdots & \ddots & \vdots\\
#7 & #8 & \dots & #9\end{array}
\right)}
\newcommand{\matricelongue}[9]	%matrice longue$
{\left(\begin{array}{ccccc} #1 & #2 &\dots & \dots &  #3\\#4
& #5
& \dots & \dots & #6\\ \vdots & \vdots  & & & \vdots\\
#7 & #8 & \dots & \dots & #9\end{array}
\right)}
\newcommand{\matricehaute}[9]	%matrice haute$
{\left(\begin{array}{cccc} #1 & #2 &\dots &  #3\\#4
& #5
& \dots & #6\\ \vdots & \vdots &  & \vdots\\\vdots & \vdots &  &
\vdots\\
#7 & #8 & \dots & #9\end{array}
\right)}

\begin{document}
%Titre
\title{Convergence rate for the method of moments with linear closure relations}
\author{Y. Bourgault${}^1$, D. Broizat${}^2$, P.-E. Jabin${}^3$}
\footnotetext[1]{Department of Mathematics and Statistics, University of Ottawa, 585 King Edward Ave, Ottawa, Ontario, Canada. E-mail: 
\texttt{ybourg@uottawa.ca}}

\footnotetext[2]{Laboratoire J.-A. Dieudonn\'e, Universit\'e de Nice --
  Sophia Antipolis, Parc Valrose, 06108 Nice Cedex 02, France, E-mail:
  \texttt{broizat@unice.fr}}

\footnotetext[3]{CSCAMM and Department of Mathematics, University of
  Maryland, College Park, MD 20742-4015, USA,
E-mail: \texttt{pjabin@cscamm.umd.edu}}
\date{}

\maketitle

\begin{abstract}We study linear closure relations for the moments' method applied to simple kinetic 
equations. The equations are linear collisional models (velocity jump processes) which are well suited to this type of approximation. In this simplified, 1 dimensional setting, we are able to  
prove stability estimates for the method (with a kinetic 
interpretation by a BGK model). Moreover we are also  
able to obtain convergence rates which automatically increase with the smoothness of the initial data.
\end{abstract}
%%%%%%%%%%%%%%%%%%%%%%%%%%%%%%%%%%%%%%%%%%%%%%%%%%%%%%%%%%%%%%%%%%%%%%%%%%%%%%%%
%%%%%%%%%%%%%%%%%%%%%%%%%%%%%%%%%%%%%%%%%%%%%%%%%%%%%%%%%%%%%%%%%%%%%%%%%%%%%%%%
%%%%%%%%%%%%%%%%%%%%%%%%%%%%%%%%%%%%%%%%%%%%%%%%%%%%%%%%%%%%%%%%%%%%%%%%%%%%%%%%
\section{Introduction}

\subsection{Quick presentation of the moments' methods for kinetic equations}
We study here an unusual but very simple choice of closure for the moments' method where we can completely characterize the stability and convergence rates of the approximation. As far as we know this very simple situation was never considered before. Before presenting this choice though, let us very briefly give the main ideas behind the method of moments.

Moments' methods have been introduced in \cite{Gr} in the context of the Boltzmann equation. This well known equation is posed on the density $f(t,x,v)$ of particles in the phase space and reads
\begin{equation}
\partial_t f+v\cdot\nabla_x f=Q(f),\label{boltz}
\end{equation}
where $Q$ is a non linear operator (in the velocity variable $v$) which expresses how the velocity of a particle may change when it has a random collision. 

Solving numerically an equation like \eqref{boltz} is in general very costly. The structure of the righthand side $Q$ is non local in $v$ and moreover the equation is posed in phase space which means that one has to work in dimension $2d+1$ if $x,v\in \R^d$ ($7$ for instance if $x\in \R^3$).

The moments' method is one possible answer to this problem and it consists in solving equations on the polynomial moments of $f$. Let $m(v)$ be polynomial in $v$ then
\begin{equation}
\partial_t <m,f>+\nabla_x\cdot<v\,m(v),f>=<m,Q(f)>,\label{momeq}
\end{equation}
 where we denote
\begin{equation}
<m,f>=\int m(v)\,f(t,x,v)\,dv.\label{defmom}
\end{equation}
Now instead of solving one equation in dimension $2d+1$, one has to solve several equations but in dimension $d+1$. Moreover in general one can expect that $<m,\,Q(f)>$ is not too complicated to compute.

However the system given by \eqref{momeq} is not closed as $vm$ is always one degree higher than $m$. Therefore no matter how many moments $m_i$ one chooses, it is never possible to express all the $v\,m_j$ in terms of the $m_i$. This is the closure problem and it means that \eqref{boltz} is never equivalent to \eqref{momeq} for any finite number of moments.

Instead one typically chooses a closure equation, {\em i.e.} a relation between $<vm_i,\,f>$ and the $<m_j,f>$ for those $i$ where $vm_i$ cannot be expressed in terms of the $m_j$. 

The first big difficulty for this type of method is how to choose the closure in order to ensure that the corresponding moments' system has good properties and gives a good approximation of Eq. \eqref{boltz}. This problem was of course recognized early on, see for instance \cite{Bo}, as well as the role of entropy, see \cite{Pe} among many others.

One of the first  systematic ways of finding a closure was introduced in \cite{Lev1} and \cite{LM}. It is still not easy to actually compute the relation which means that it is often computed numerically instead (see \cite{ST}, \cite{CL}, \cite{CLW}). Different closures can of course be used (see \cite{Tor2} for example).  

Theoretically even checking that the corresponding method leads to a hyperbolic system is not easy (we refer for instance to \cite{Br}, \cite{CFL}). Proving convergence rates seems to be out of reach for the time being although in practice it seems to be a good approximation (see \cite{LP} for a numerical study). 

Let us also mention that the methods of moments has also been used for theoretical purposes (as in \cite{De}) and not only numerical computations.

We conclude this very brief overview by refering to \cite{Bi}, \cite{St} or \cite{To} for more on numerical simulations for kinetic equations in this context.

%%%%%%%%%%%%%%%%%%%%%%%%%%%%%%%%%%%%%%%%%%%%%%%%%%%%%%%%%%%%%%%
\subsection{Linear closure relations}
%%%%%%%%%%%%%%%%%%%%%%%%%%%%%%%%%%%%%%%%%%%%%%%%%%%%%%%%%%%%%%%
The guiding question in this article is whether it can make sense to consider a linear closure relation. This is certainly delicate in the nonlinear case of Boltzmann eq. \eqref{boltz}. Instead we choose a simplified 1d setting where it is possible to fully analyze the method.

Instead of \eqref{boltz}, we consider the linear model 
\begin{equation}\label{eq_neutronics}
\left\{\begin{array}{lll}
\del_t f+v\del_x f=L(f),\qquad (x,v)\in\r\times \r\\\\
\ds L(f)=\int_{\r}Q(v,v^*)f(t,x,v^*)dv^*-\lambda f
\end{array}\right., 
\end{equation}
with $\lambda>0$ and where the operator $Q$ corresponds to a velocity jump process.  While much simplified with respect to \eqref{boltz}, it is not uninteresting in itself, with applications to physics (see \cite{Ch}, \cite{Win}) or biology (see \cite{ODA} for example). The equation had to be supplemented with some initial data, which for simplicity we assume to be compactly supported in velocity
\begin{equation}
f(t=0,x,v)=f^0(x,v)\in L^2(\R^2),\quad \mbox{supp}\,f^0\subset \r\times I.\label{initialdata}
\end{equation} 

In general $Q$ could even be assumed to depend on $t$ and $x$. Here we make the additional approximation
\begin{equation}\label{forme_noyau}
Q(v,v^*)=\left(q(v)\sum_{j=0}^d \alpha_j {v^*}^j\right)\1_{\{(v,v^{*})\in I^{2}\}},
\end{equation}
with $q$ smooth and compactly supported in some interval $I$, $d\in\N^*$ and $(\alpha_j)_{0\leq j\leq d}\in
\r^{d+1}$. With this special form, one of course expects to be in a very favorable situation for the method of moments. Hence this should be seen as a simple toy model where the method can easily be tested.

Denoting the moments of the solution $f$ by
\begin{equation}\label{def_moments_f}
\mu_i^f(t,x):=\int_{I}v^if(t,x,v)dv,\qquad i\in\N,
\end{equation}
Eq. \eqref{eq_neutronics} simply becomes
\begin{equation}\label{eq_neutronics_bis}
\begin{array}{llll}
\del_t f+v\del_x f= 
\ds L(f)  =  q(v)\sum_{j=0}^d \alpha_j \mu_j^f -\lambda f.
\end{array}
\end{equation}
As we work in dimension $1$, the structure of the hierarchy of equations on the moments is also very simple 
\begin{equation}\label{eq_moments}
 \del_t \mu_i^f+\del_x \mu_{i+1}^f=\gamma_i\left(\sum_{j=0}^d \alpha_j
\mu_j^f \right)-\lambda \mu_i^f,\qquad\qquad  i\in\N,
\end{equation}
where we truncate at order $N$ and we define the moments of $q$
\begin{equation}\label{def_gamma}
\gamma_i=\mu_i^q=\int_{I}v^i q(v)dv,\qquad i\in\N. 
\end{equation}
In order to close the system, it would be necessary to be able to express $\mu_{N+1}^f$ in terms of the $\mu_i^f$ for $i\leq N$.  
The linear closure relation that we study here consists in assuming that $\mu_{N+1}^f$ is a linear combination of the lower moments.

That means that instead of \eqref{eq_neutronics_bis} or \eqref{eq_moments}, we solve
\begin{equation}\label{method_moments}
\left\{\begin{array}{lll}
   \del_t \mu_i+\del_x \mu_{i+1}=\ds\gamma_i\left(\sum_{j=0}^d
\alpha_j \mu_j\right)-\lambda \mu_i,\qquad i=0,\dots,N\\
\mu_{N+1}=\ds\sum_{i=0}^N a_i \mu_i.
  \end{array}\right.,
\end{equation}
with $(a_0,\dots,a_N)$ given real coefficients that one should choose in the ``best'' possible way.

One can rewrite $\eqref{method_moments}$ in matrix form as
\begin{equation}\label{systeme_hyperbolique}
 \del_t M_N+\del_x (A M_N)=BM_N,
\end{equation}
where $M_N=M_N(t,x)=(\mu_0(t,x),\dots,\mu_N(t,x))^{T}\in \r^{N+1},$ and $A,B\in
M_{N+1}(\r)$ are defined by
\begin{equation}\label{matrices_systeme}
A=\left(\begin{array}{cccccc}
0 & 1   &  & \\ \vspace{2mm}
 &   \ddots & \ddots & \\\vspace{2mm}
 &    & 0 & 1\\\vspace{2mm}
a_0   & \dots & a_{N-1} & a_N
\end{array}
\right),\quad
B=\left(\begin{array}{ccc|ccc}
 &  \vdots  &  & & 0 & \\ 
 &    &  & & & \\ 
 &  (\gamma_i \alpha_j)_{\genfrac{}{}{0pt}{}{0\leq i\leq N}{0\leq j\leq d}}  & 
& & \vdots &\\ 
&    &  & &  &\\
& \vdots   &  & & 0 &
\end{array}
\right)-\lambda I_{N+1}.
\end{equation}
%%%%%%%%%%%%%%%%%%%%%%%%%%%%%%%%%%%%%%%%%%%%
%%%%%%%%%%%%%%%%%%%%%%%%%%%%%%%%%%%%%%%
\subsection{Basic properties of the method and the main result}
%%%%%%%%%%%%%%%%%%%%%%%%%%%%%%%%%%%%%%%
%%%%%%%%%%%%%%%%%%%%%%%%%%%%%%%%%%%%%%%%%%%%
Given the simple form of \eqref{method_moments} or \eqref{systeme_hyperbolique}, some properties are obvious.  
Given its linear nature, the system is hyperbolic if the characteristic polynomial
\begin{equation}\label{poly_carac_A}
\chi_A(X)=det(XI_{N+1}-A)=X^{N+1}-\sum_{i=0}^N a_i X^i 
\end{equation}
has $N+1$ real roots, denoted by
\begin{equation}\label{spec_A}
spec(A)=\{\lambda_0,\dots,\lambda_N\} .
\end{equation}
This is enough to guarantee the well posedness of the numerical system \eqref{systeme_hyperbolique} but not necessarily good stability properties. Norms of the numerical appoximation could for instance grows fast as $N$ increases. If the initial data is compactly supported in $x$ then the solution is as well and the support propagates with speed $\max_k |\lambda_k|$.

On the other hand, the major inconvenients of such a method are also pretty clear. For instance positivity of the even moments will likely not be preserved. Even the positivity of the macroscopic density $\mu^f_0$ has no reason to be propagated. 

However a careful analysis can in fact reveal that it is possible to choose appropriately the coefficients $a_i$ in order to have not only stability but also very fast convergence.

For simplicity, assume that $I=[-1,\ 1]$ (just rescale and translate otherwise) and choose the Tchebychev polynomial for $\chi_A$ or
\begin{equation}
\forall 0\leq k\leq N,\qquad \lambda_{k}=\cos\left(\left(\frac{2k+1}{2N+2}\right)\pi\right),
\label{tcheblambda}
\end{equation}
or 
\begin{equation}
\chi_{A}^{(N+1)}(X)=\prod_{k=0}^{N}\left(X-\cos\left(\left(\frac{2k+1}{2N+2}\right)\pi\right)\right).\label{tchebchi}
\end{equation}
Then it is possible to show:
\begin{Th}\label{th stabilite tchebychev} 
Assume that $f^0$ satisfies \eqref{initialdata}, that $Q$ satisfies \eqref{forme_noyau}. Then the solution to the truncated moments' hierarchy \eqref{method_moments} where the coefficients $a_i$ are chosen according to \eqref{tcheblambda} or \eqref{tchebchi} satisfies
\begin{equation}\label{estimation_de_stabilite_macro_tchebychev}
\sup_{0\leq i\leq N}\sup_{t\in[0,T]}\|\mu_{i}(t)\|_{L^{2}(\r)}\leq 
 e^{T C_{d,\alpha}\|q\|_{L^{2}(I)}}\|f^{0}\|_{L^{2}(\r^{2})},
\end{equation} 
where $C_{d,\alpha}=\sqrt{\pi(d+1)}\left(\sum_{j=0}^{d}{\alpha_{j}}^{2}\right)^{1/2}$ is independent of $N$.

In addition if $f^0\in H^k(\r^{2})$, $d=0$, denoting by $f$ the corresponding solution to \eqref{eq_neutronics} and defining its moments by \eqref{def_moments_f}, one has
\begin{equation}\label{errorestimatebasic}
\|\mu_{0}-\mu_{0}^{f}\|_{L^{\infty}\left([0,T],L^{2}(\r)\right)}\leq \frac{C}{N^{k-3/4}}\times\|f^{0}\|_{H^{k}(\r^{2})}.	
\end{equation}
where $C\geq 0$ depends on $T,\lambda,q$ and $k$.
\end{Th}
\begin{Rk} The convergence result is given for the first moment and $d=0$. Higher moments would lose small powers of $N$ so the same proof would give
\[
\|\mu_{i}-\mu_{i}^{f}\|_{L^{\infty}\left([0,T],L^{2}(\r)\right)}\leq \frac{C\,N^{k\,(i+d)/N}}{N^{k-1}}\times\|f^{0}\|_{H^{k}(\r^{2})}.	
\]
We do not know whether those results are optimal or not and in particular whether the numerical value $3/4$ in \eqref{errorestimatebasic} is optimal (though actual numerical simulations do suggest it is).
\end{Rk}
Hence in this setting, the conclusion of this analysis is that the linear method of moments should be seen in the same light as spectral methods (see \cite{ShT} for instance). It is stable, converges automatically at order $k-3/4$ if the initial data is in $H^k$ but does not propagate any additional property (positivity being probably the most important).

The next section is a more detailed presentation of the stability and convergence analysis in the general case (without necessarily choosing the Tchebychev points). The corresponding technical justifications and calculations are presented in two separate sections. Theorem \ref{th stabilite tchebychev} is proved in section 5. The last section is an appendix and recalls some well known results.
%%%%%%%%%%%%%%%%%%%%%%%%%%%%%%%%%%%%%%%%%%%%%%%%%%%%%%%%%%%
%%%%%%%%%%%%%%%%%%%%%%%%%%%%%%%%%%%%%%%%%%%%%%%%%%%%%%
\section{Stability and convergence results}
%%%%%%%%%%%%%%%%%%%%%%%%%%%%%%%%%%%%%%%%%%%%%%%%%%%%%%%%%
%%%%%%%%%%%%%%%%%%%%%%%%%%%%%%%%%%%%%%%%%%%%%%%%%%%%%%%%%%%%%
We present here in more details the kind of stability and convergence results that can be proved for the method \eqref{method_moments} without assuming any particular choice of the eigenvalues (like \eqref{tcheblambda}). We give here the main ideas in the approach and leave the technical proofs to further sections.
%%%%%%%%%%%%%%%%%%%%%%%%%%%%%%%%%%%%%%%%%%%%%%%%%%
%%%%%%%%%%%%%%%%%%%%%%%%%%%%%%%%%%%%%%%%%%%%%%%%%%%%
\subsection{Eigenvectors for System \eqref{systeme_hyperbolique}}
%%%%%%%%%%%%%%%%%%%%%%%%%%%%%%%%%%%%%%%%%%%%%%%%%%%%%%%%%
%%%%%%%%%%%%%%%%%%%%%%%%%%%%%%%%%%%%%%%%%%%%%%%%%%%%%%%%%%%
As they enter in some of the estimates, we start by a short parenthesis about the eigenvectors for the matrix $A$.
Define $P\in M_{N+1}(\r)$
the matrix of the
eigenvectors; then 
\begin{equation}\label{matrice_de_passage}
P_{i,j}=\lambda_{j}^i,\qquad 0\leq i,j\leq N.
\end{equation}
Its inverse $P^{-1}$ can be computed
as easily and
\begin{equation}\label{matrice_de_passage_inverse}
P_{i,j}^{-1}=\frac{
\tilde p_{i,j}}{\pi_i},\qquad 0\leq i,j\leq N,
\end{equation}
where
\begin{equation}\label{formules_P_-1}
\pi_i=\prod_{j\neq i}
(\lambda_i-\lambda_j),\qquad
\tilde p_{i,j}=(-1)^{N-j}
\sum_{\genfrac{}{}{0pt}{}{k_1<\ldots<k_{N-j}}{k_l\neq i\ \forall l}} 
\prod_{l=1}^{N-j} \lambda_{k_l}=\sum_{l=0}^j \frac{a_l}{\lambda_i^{j-l+1}}.
\end{equation}
with the convention that $\tilde p_{i,N}=1$.

Moreover, an easy computation shows that:
\begin{equation}\label{autre formule inverse P}
\tilde{p}_{i,j}=\lambda_i^{N-j}-a_N\lambda_i^{N-j-1}-\dots-a_{j+2}\lambda_i-a_{j+1},\qquad 0\leq i,j\leq N.
\end{equation}
We can notice that $\tilde{p}_{i,j}$ is an homogeneous polynomial of degree $N-j$ in the eigenvalues 
$(\lambda_0,\dots,\lambda_N)$.
%
%
%%%%%%%%%%%%%%%%%%%%%%%%%%%%%%%%%%%%%%%%%%%%%%%%%%%%%%%%%%%%%%%
%%%%%%%%%%%%%%%%%%%%%%%%%%%%%%%%%%%%%%%%%%%%%%%%%%%%%%%%%%%%%%
\subsection{Stability estimate and kinetic interpretation of the method}
%%%%%%%%%%%%%%%%%%%%%%%%%%%%%%%%%%%%%%%%%%%%%%%%%%%%%%%%%%%%%
%%%%%%%%%%%%%%%%%%%%%%%%%%%%%%%%%%%%%%%%%%%%%%%%%%%%%%%%%%%%%%%
Let us first state our main stability estimate.
\begin{Th}\label{th stabilite} 
Assume \eqref{initialdata} and that $\{\lambda_{0},\dots,\lambda_{N}\}\subset I$. Moreover assume that there exists
a function $\rho_N(v)$, positive on $I$, such that
\begin{equation}\label{poids_rho_N}
\int_{I}\frac{R(v)}{\rho_{N}(v)}dv=\sum_{k=0}^N R(\lambda_k),\qquad \forall R\in\r_{2N+1}[X].
\end{equation}
Then, the hyperbolic system \eqref{systeme_hyperbolique} is stable and
\begin{equation}\label{estimation_de_stabilite_macro}
\sup_{t\in[0,T]}\|\mu_{i}(t)\|_{L^{2}(\r)}\leq 
\left(\sum_{k=0}^{N}\lambda_{k}^{2i}\right)^{1/2}\times e^{T C_{N}(q)}
\times C_{N}(f^{0}),
\qquad i=0,\dots,N	
\end{equation}
where

\begin{equation}\label{constante_pour_donnee_initiale}
C_{N}(f^{0})=\left(\int\!\!\!\int_{\r\times I} |f^{0}(x,v)|^{2}\rho_{N}(v)dxdv\right)^{1/2},
\end{equation}

\begin{equation}\label{constante_pour_q}
C_{N}(q)=\Lambda_{N,d}\left(\int_{I}|q(v)|^{2}\rho_{N}(v)dv\right)^{1/2}-\lambda,
\end{equation}

\begin{equation}\label{constante Lambda}
\Lambda_{N,d}=\sqrt{\sum_{j=0}^{d}\alpha_{j}^{2}}\sqrt{\sum_{j=0}^{d}\sum_{k=0}^{N}\lambda_{k}^{2j}}
=\sqrt{\sum_{j=0}^{d}\alpha_{j}^{2}}\sqrt{\sum_{k=0}^{N}\frac{1-\lambda_{k}^{2d+2}}{1-\lambda_{k}^{2}}}.	
\end{equation}
\end{Th}
This result does not assume any particular distribution on the eigenvalues but of course it is by no means guaranteed in general that one could find $\rho_N$ satisfying \eqref{poids_rho_N}. Notice that the corresponding relation is really a quadrature formula for computing integrals on $I$ which we ask to be exact for polynomials of degree $2N+1$.

We do not prove this result directly on the system \eqref{systeme_hyperbolique}. Instead we show a corresponding result on a linear BGK problem (see \cite{BGK} for the simplification of collision kernels into BGK models).

For any $\eps>0$ and $N\geq d$, consider the equation 
\begin{equation}\label{BGK}
\begin{split}
&\del_t f_N^{\eps}+v\del_x{f_N^\eps}  =  \frac{M^{(N)}\,f_N^{\eps}-f_N^{\eps}}{\eps}
+L(f_N^{\eps})\\
&f_N^{\eps}(0,x,v)=f_0(x,v), 
\end{split}\end{equation}
where the linear operator $L$ is defined by $\eqref{eq_neutronics_bis}$
and the Maxwellian $M^{(N)}:f\mapsto  M^{(N)}f$ satisfies the moment conditions
\begin{equation}\label{conditions moments maxwellienne}
\left\{\begin{array}{llll}
\ds\int_{I}v^i M^{(N)}\,f dv & = & \ds\int_{I}v^i f dv &  i=0,\dots,N,\\\\
\ds\int_{I}v^{N+1} M^{(N)}\,f dv & = &\ds \sum_{i=0}^N a_i \int_{I}v^i f dv & 
\end{array}
\right..
\end{equation}
This problem is a kinetic approximation of the macroscopic problem
\eqref{method_moments} as formally when $\eps\longrightarrow 0$ then one recovers \eqref{method_moments} from \eqref{BGK}. 

Indeed,
multiplying \eqref{BGK} by $v^i$ and integrating over
$v\in I$, we obtain 
\[
\del_t \mu_i^{\eps}+\del_x\mu_{i+1}^{\eps}=
\gamma_i\left(\sum_{j=0}^{d} \alpha_j \mu_j^{\eps}\right)-\lambda \mu_i^{\eps}
\quad \mbox{for}\ i=0,\dots,N,\]
where $\mu_i^{\eps}:=\int_{I}v^i f_N^{\eps} dv$ for
$i\in\N$.
Moreover, when $\eps\rightarrow 0$, we have
formally
\begin{equation*}
\mu_{N+1}^{\eps}=\int_{I}v^{N+1}f_N^{\eps} dv\sim
\int_{I}{v^{N+1}}M^{(N)} f_N^{\eps} dv=\sum_{i=0}^N a_i \int_{I}v^i f_N^{\eps}
dv=\sum_{i=0}^N a_i \mu_i^{\eps},
\end{equation*}
which is our closure relation. 

The interest of \eqref{BGK} is to make some computations more transparent and easier to follow and in addition at the kinetic level which is more natural given the original equation. If we can prove stability estimates for \eqref{BGK} that are uniform in $\eps$ then simply by passing to the limit, we will obtain estimates for \eqref{systeme_hyperbolique}-\eqref{method_moments}.

%The aim is to find a norm which provides bounds on the solutions
%$(f_N^{\eps}(t))_{\eps>0}$ in function of $f^0$, uniformly in
%$\eps$. It will give a stability estimate for the problem
%$\eqref{method_moments}$ at kinetic level. 

The
most obvious choice for the maxwellian is simply

\begin{equation}\label{formule_maxwellienne}
M^{(N)}\,f=\sum_{i=0}^N \left(\int_{I}v^i f dv\right)m_i, 
\end{equation}
where the $v\mapsto m_i(v)$, $i=0,\dots,N$, are any functions satisfying 
\begin{equation}\label{cond_moments_m_i}
\int_{I}v^j m_i(v) dv=\delta_{i,j},\quad 0\leq j\leq N,\qquad\qquad
\int_{I}v^{N+1}m_i(v) dv =a_i. 
\end{equation}
The conditions \eqref{cond_moments_m_i} ensure
that \eqref{conditions moments maxwellienne} holds.

There are obviously many ways to choose the $m_i$ s.t. \eqref{cond_moments_m_i} is satisfied. What we are looking for is a choice compatible with an inner product $\Phi_N$ such that the application $M^{(N)}$ is an orthogonal projection for $\Phi_N$. Formally this implies that for any $f$
\[
\Phi_N(f,M^{(N)} f)\leq \Phi_N(f,f).
\]
In addition this inner product should have good symmetry properties s.t. for  $f$ and $g$
\[
\Phi_N(f,v g)=\Phi_N(vf,g).
\]
The simplest way to ensure this is to look for a weight $\rho_N$ s.t.
\[
\Phi_N(f,g)=\int_I f(v)\,g(v)\,\rho_N(v)\,dv.
\]
In that case formally
\[
\Phi_N(f,v\partial_x f)=\frac{1}{2}\partial_x\,\Phi_N(f,vf).
\]
Therefore if $f_N^\eps$ solves \eqref{BGK} then one expects that
\[
\frac{d}{dt}\int_\R \Phi_N(f_N^\eps,f_N^\eps)\,dx\leq 2\,\int_\R \Phi_N(f_N^\eps,\,L(f_N^\eps))\,dx.
\]
This is the strategy that we implement. Find appropriate conditions on $\rho_N$ and the $m_i$ to obtain the correct structure and then simply bound 
$\Phi_N(f, L(f))$ in terms of $\Phi_N(f)$. 

%Since the map $M$ acts only on the velocity variable of $f=f(t,x,v)$, we can see
%$M$ as a linear map 
%
%\begin{equation}\label{M_endo_E}
%M:E\longrightarrow E,
%\end{equation}
%
%where $E$ is a Banach space of functions of the
%velocity variable $v$.

%A formal computation shows the interest of having $\|Mf\|_E\leq
%\|f\|_E$ (i.e. $M$ being a shrinking operator). Indeed, if we omit the transport operator $v\del_x f$, the
%interaction term $L(f)$ and simply study the
%equation $\del_t f =\frac{1}{\eps}(Mf-f)$, then we have
%
%\[\begin{split}
%\|f(t+h)\|_E&\underset{h\rightarrow
%0^+}{\sim}\left\|f(t)+\frac{h}{\eps}(Mf(t)-f(t))\right\|_E\\
%&\leq
%\left(1-\frac{h}{\eps}\right)\left\|f(t)\right\|_E+\frac{h}{\eps}\|Mf(t)\|_E.
%\end{split}
%\]
%
%Thus $$\frac{\|f(t+h)\|_E-\|f(t)\|_E}{h}\leq
%\frac{\|Mf(t)\|_E-\|f(t)\|_E}{\eps}\leq 0,$$
%
%Hence, letting $h\rightarrow 0^+$, $\|Mf\|_E\leq\|f\|_E$ formally implies
%$\del_t
%\|f(t)\|_E\leq 0$, i.e. the equation operator is formally dissipative, which
%would
%provide the wanted stability bounds on $\|f(t)\|_E$.

%Therefore, the aim is now to build a Banach space $E$ and a Maxwellian $M$ such that 
%
%\begin{equation}
%\|Mf\|_E \leq \|f\|_E,\qquad f\in E
%\end{equation}
%
%and such that $\|.\|_E$ is compatible with the transport operator
%$T=\del_t+v\del_x$. 

%In fact, we can build an Hilbert space $E$ with these properties
For the first part of this strategy we actually prove
\begin{Th}\label{th maxwellienne} Let $\rho_{N}(v)$ a positive function on $I$ such that
\begin{equation}\label{conditions moments poids}
\int_{I}\frac{R(v)}{\rho_{N}(v)}dv=\sum_{k=0}^N R(\lambda_k),\qquad \forall R\in\r_{2N+1}[X].
\end{equation}
We set
\begin{equation*}
E_{N}=L^2(I,\rho_{N}(v)dv)=\left\{f\,\, measurable,\quad \int_{I}|f(v)|^2\rho_{N}(v)dv<\infty\right\}.
\end{equation*}
Then, the map $\phi_{N}$ defined by
\begin{equation}\label{weightedl2}
\phi_{N}(f,g)=\int_{I}f(v)g(v)\rho_{N}(v)dv,\qquad (f,g)\in E_{N}^2
\end{equation}
is an inner product on $E_{N}$, and the
Maxwellian $M^{(N)}:E_{N}\rightarrow E_{N}$, defined by
$$M^{(N)}f=\sum_{i=0}^N\left(\int_{I}v^{i}f(v)dv\right)\frac{\tilde{T_i}(v)}{\rho_{N}(v)},\qquad f\in E_{N},$$
is an orthogonal projection %(and thus a shrinking operator for the hilbertian norm on $E_{N}$), 
and satifies \eqref{conditions moments maxwellienne},
where $(\tilde{T_i})_{0\leq i\leq N}$ is the basis of $\r_N[X]$ defined by:
$$\tilde{T_i}(X)=\sum_{k=0}^NQ_{k,i}X^k,\qquad 0\leq i\leq N.$$
with $Q=(P^{T})^{-1}P^{-1}$.

Furthermore, with that choice of $M^{(N)}$, any solution $f_{N}^{\eps}$ of the problem \eqref{BGK} formally satisfies:
\begin{equation}\label{bornes_stab_BGK}
\frac{d}{dt}\int_{\r}\phi_{N}(f_N^{\eps},f_N^{\eps})dx\leq 2
\int_{\r}\phi_{N}(f_N^{\eps},L(f_N^{\eps}))dx. 
\end{equation}
\end{Th}
%%%%%%%%%%%%%%%%%%%%%%%%%%%%%%%%%%%%%%%%%%%%%%%%%%%%%%%%%%%%%%%%
%%%%%%%%%%%%%%%%%%%%%%%%%%%%%%%%%%%%%%%%%%%%%%%%%%%%%%%%%%%%%%%%
\subsection{Error estimate}
%%%%%%%%%%%%%%%%%%%%%%%%%%%%%%%%%%%%%%%%%%%%%%%%%%%%%%%%%%%%%%%%%
%%%%%%%%%%%%%%%%%%%%%%%%%%%%%%%%%%%%%%%%%%%%%%%%%%%%%%%%%%%%%%%%%
We now turn to the convergence of the method of moments to the solution. For simplicity we state the results here for the case $d=0$ in \eqref{forme_noyau}, namely 
$$Q\left(v,v^{*}\right)=q(v)\1_{\{(v,v^{*})\in I^{2}\}}.$$ 
We also assume that $I\subset [-1,\ 1]$, still for simplicity.

The convergence results of course require a stability estimate. However in themselves, they do not use the specific form of the closure. 

Therefore here we do not assume any specific closure relation. Instead we assume that we have a well defined methods of moments, {\em i.e.}, some way of computing $\mu_i$ which satisfy  
\begin{equation}\label{method_moments_d=0}
\begin{array}{lll}
   \del_t \mu_i+\del_x \mu_{i+1}=\ds\gamma_i
\mu_0-\lambda \mu_i,\qquad i=0,\dots,N\\
\mu_{N+1}=F(\mu_1,...,\mu_N).
  \end{array},
\end{equation}
where
$$\gamma_{i}=\int_{I}v^{i}q(v)dv.$$
Moreover we assume that the corresponding method has good stability estimates in the sense that 
\begin{equation}\label{estimation_de_stabilite_macro_d=0_sous_hypotheses}
\begin{split}
& \exists C,\,\gamma\geq 0,\ \mbox{For any}\ (\mu_i)_{i=0..N}\ \mbox{and}\ (\tilde\mu_i)_{i=0..N} \mbox{ solutions to \eqref{method_moments} for}\ f^0,\;\tilde f^0,\\
&  \mbox{then}\qquad \sum_{i=0}^{N}\|\mu_{i}-\tilde\mu_i\|_{L^{\infty}([0,T],\;H^k(\r))}\leq 
CN^{\gamma}\|f^{0}-\tilde f^0\|_{L^2(I,\ H^k(\r))}, 
\end{split}\end{equation}
where $C,\;\gamma$ may depend on $T$, $q$, $\lambda$ but not on $N$ or $f^0$.

Note that Th. \ref{th stabilite} can indeed be expected to imply such a result. The exponent $\gamma$ will depend on the choice of the coefficients in our method. Of course Th. \ref{th stabilite} controls only the $L^2$ norm of the $\mu_i$. However as the method \eqref{method_moments} is linear, the $\partial_x^k \mu_i$ are also a solution to the same system and an estimate like \eqref{estimation_de_stabilite_macro_d=0_sous_hypotheses} can be derived. For a more detailed analysis of how to obtain \eqref{estimation_de_stabilite_macro_d=0_sous_hypotheses} from Th. \ref{th stabilite}, we refer to Section \ref{Tcheby} where it is performed when the $\lambda_k$ are the Tchebychev points.

For any method that satisfies \eqref{estimation_de_stabilite_macro_d=0_sous_hypotheses}, then one has the following convergence result
\begin{Th} Assume \eqref{initialdata} with $I\subset [-1,\ 1]$, that the method \eqref{method_moments_d=0} satisfies \eqref{estimation_de_stabilite_macro_d=0_sous_hypotheses} for some $\gamma\geq 0$ and that $f^{0}\in H^{k}(\r^{2})$, with $k\in\N^{*}$.  Consider $f$ the solution to \eqref{eq_neutronics} with $d=0$ and the corresponding solution $\mu_i$ to \eqref{method_moments_d=0}. Then, for all $T\geq 0$ and for all $N\geq 1$, 
we have the estimate
\begin{equation}\label{estimation_erreur}
\|\mu_{0}-\mu_{0}^{f}\|_{L^{\infty}\left([0,T],L^{2}(\r)\right)}\leq \frac{C}{N^{k-\gamma}}\times\|f^{0}\|_{H^{k}(\r^{2})}.	
\end{equation}	
where $C\geq 0$ depends on $T,\lambda,q$ and $k$ but not on $N$.
\label{errorestimate}
\end{Th}	
This error estimate is a sort of interpolation between the stability bounds and the following result for $C^\infty$ solution to \eqref{method_moments_d=0}
\begin{Prop}\label{resultat_smooth_case} Assume that the $(\mu_i)_{i=0...N+1}$ solve
\begin{equation}
\partial_t \mu_i+\partial_x \mu_{i+1}=\gamma_i\,\mu_0-\lambda\mu_i,\qquad i=0,..., N,
\end{equation}
 with $\mu_i(t=0)=0$ for $i=0...N$.
Then, for all $T\geq 0$ and for all $N\geq 1$, 
we have the estimate
\begin{equation}\label{formule erreur cas_regulier}
\|\mu_0-\mu_0^{f}\|_{L^{\infty}([0,T],L^2(\r))}\leq C\times 
\frac{T^N}{N!}\,\sum_{i=0}^N
\|\del_x^N \mu_{i}\|_{L^\infty([0,T],\ L^2(\r))},
\end{equation}
where $C\geq 0$ depends on $T,\,\lambda,\,q$.
\end{Prop}
%%%%%%%%%%%%%%%%%%%%%%%%%%%%%%%%%%%%%%%%%%%%%%%%%%%%%%%%%%%%%%%%%%%%%%%%%%%%%%%%%%%%%%%%%%%%%%%%%%%%%%%%%%%%%%%%%%%%%%%%%%%%%%%%%%%%%%%%%%%%%%%%%%%%%%%%%%%%%%%
\section{Proof of Theorems \ref{th stabilite} and \ref{th maxwellienne}}
%%%%%%%%%%%%%%%%%%%%%%%%%%%%%%%%%%%%%%%%%%%%%%%%%%%%%%%%%%%%%%%%%%%%%%%%%%
%%%%%%%%%%%%%%%%%%%%%%%%%%%%%%%%%%%%%%%%%%%%%%%%%%%%%%%%%%%%%%%%%%%%%%%%%%%
To simplify the presentation, we omit here the subscript $N$ in
$E_{N}$, $\phi_{N}$, $\rho_{N}$, and the superscript $N$ in $M^{(N)}$. 
%%%%%%%%%%%%%%%%%%%%%%%%%%%%%%%%%%%%%%%%%%%%%%%%%%%%%%%%%%%%%%%5
\subsection{Elementary space decomposition}
%%%%%%%%%%%%%%%%%%%%%%%%%%%%%%%%%%%%%%%%%%%%%%%%%%%%%%%%%%%%%%%%%%%
The difficulty is to combine the fact that $M$ has to be an orthogonal projection for $\Phi$ with the symmetry property on $\Phi$. We take here a slighty more general approach by not assuming directly that $\Phi$ satisfies \eqref{weightedl2}.

Consider a Maxwellian which has the form
\eqref{formule_maxwellienne} and satisfies \eqref{cond_moments_m_i}. 
First,  such an application $M$ is a projection
because $M\circ M=M$, which is a straightforward consequence of
\begin{equation*}%label{points_fixes_M}
M(m_i)=\sum_{j=0}^N \left(\int_{I}v^j m_i dv\right)m_j=m_i,\qquad
i=0,\dots,N.
\end{equation*}
Moreover, one has  that 
\begin{equation*}%\label{noyau_image_M}
\begin{array}[]{llll}
Ker M=\ds \left\{f\in E,\quad \int_{I}v^if(v)dv=0, \quad
i=0,\dots,N \right\}:=K,\\\\
Ker(M-I)=Im(M)=\ds Span(m_0,m_1,\dots,m_N):=V,
\end{array}
\end{equation*}
and we have the space decomposition 
\begin{equation*}%\label{somme directe}
E=V\oplus K,
\end{equation*}
with 
\begin{equation*}%\label{dim_espaces_proj}
dim(V)=codim(K)=N+1.
\end{equation*}
We start by a more detailed explanation of the sufficient conditions to obtain Th. \ref{th maxwellienne} 
\begin{Lemma}\label{estimation_formelle} Assume that the 
inner product $\phi:E\times E\rightarrow \r$
satisfies 
\begin{equation*}%\label{cond_orthogonalite}
K=V^{\bot,\phi}, \mbox{(i.e. the decomposition $V\oplus K$ becomes orthogonal)}.
\end{equation*}
\begin{equation}\label{cond_symetrie}
\forall (f,g)\in E^2,\qquad \phi(vf,g)=\phi(vg,f).
\end{equation}
Then, for any solution $f_N^{\eps}=f_N^{\eps}(t,x,v)$ of the problem
$\eqref{BGK}$, inequality \eqref{bornes_stab_BGK} formally holds:
\begin{equation*}%\label{rappel_bornes_stab_BGK}
\frac{d}{dt}\int_{\r}\|f_N^{\eps}(t,x,.)\|^2_{\phi}dx\leq 2
\int_{\r}\phi(f_N^{\eps},L(f_N^{\eps}))dx. 
\end{equation*}
\end{Lemma}

\noindent\underline{Proof} : 

Take a smooth $f=f(t,x,v)$ such that $\del_t f+v\del_x f =\frac{1}{\eps}(Mf-f)+L(f)$. We
have
\[
 \begin{split}
&\frac{d}{dt} \int_{\r}\|f(t,x,.)\|^2_{\phi}dx  =
 2\int_{\r}\phi(f,\del_t
f)dx\\
&\ =  \ds 2\left(\frac{1}{\eps}\int_{\r}\phi(f,Mf-f)dx-\int_{\r}\phi(f,
v\del_x f)dx+\int_{\r}\phi(f,L(f))dx\right).
 \end{split}
\]
Since $Mf-f\in K$ and $Mf\in V$, we have $\phi(Mf,Mf-f)=0$, thus
\[
\phi(f,Mf-f)=\phi(f-Mf,Mf-f)=-\|f-Mf\|_{\phi}\leq0.
\]
We deduce 
\begin{equation*}%\label{bornes sur f}
\frac{d}{dt}\int_{\r}\|f(t,x)\|^2_{\phi}dx\leq
-2\int_{\r}\phi(f,v\del_x f)dx+2\int_{\r}\phi(f,L(f))dx. 
\end{equation*}
Thus, having $\int_{\r}\phi(f,v\del_x f)dx=0$ is sufficient to obtain
\eqref{bornes_stab_BGK}. 

Since $\phi(f,v\del_x f)=\del_x \phi(f,vf)-\phi(\del_x f, vf)$, we can write
\begin{equation*}%\label{integrale a annuler}
\int_{\r}\phi(f,v\del_x
f)dx=\frac{1}{2}\int_{\r}\left(\phi(f,v\del_x f)-\phi(vf,\del_x
f)\right)dx. 
\end{equation*}
Therefore the symmetry condition $\eqref{cond_symetrie}$ on $\phi$
is enough to conclude. If $f$ is not smooth enough to follow the previous steps, one simply regularizes it by convolution in $x$. As the equation is linear and $x$ is only a parameter in $M$ and $L$ then the regularized function solves the same equation. Therefore it satisfies \eqref{bornes_stab_BGK} and letting the regularizing parameter vanish, one recovers the same inequality for $f$.
\cqfd

\bigskip

Now, take a inner product $\phi$ on $E$ such that
$K=V^{\bot,\phi}$, {\em i.e.} $K$ and $V$ are orthogonal for this inner product. The symmetry condition $\eqref{cond_symetrie}$ is obviously equivalent to
\begin{equation}
\label{C1}\forall (f,g)\in V^2,\qquad
\phi(vf,g)=\phi(vg,f),
\end{equation}
\begin{equation}
\label{C2}\forall (f,g)\in K\times V,\qquad \phi(vf,g)=\phi(vg,f),
\end{equation}
\begin{equation}
\label{C3}\forall (f,g)\in K^2,\qquad \phi(vf,g)=\phi(vg,f).
\end{equation}
We study each of those conditions in the following subsections.
%%%%%%%%%%%%%%%%%%%%%%%%%%%%%%%%%%%%%%%%%%%%%%%%%%%%%%%%%%%%%%%%%%%%%%%%%%%%%%%%%%%%%%%%%%%%%%%%%%%%%%%%%%%%%%%%%%%%%%%%%%%%%%%%%%%%%%%%%%%%%%
\subsection{Study of the condition \eqref{C1}}
%%%%%%%%%%%%%%%%%%%%%%%%%%%%%%%%%%%%%%%%%%%%%%%%%%%%%%%%%%%%%%%%

\begin{Prop} The condition \eqref{C1} is equivalent to 
\begin{equation*}%\label{cond sur Q}
A^{T}Q=QA, 
\end{equation*}
where $A$ is the matrix defined by $\eqref{matrices_systeme}$ 
and $Q$ is the symmetric definite positive matrix defined by 
\begin{equation*}%\label{def Q}
Q_{i,j}=\phi(m_i,m_j),\qquad 0\leq i,j\leq N. 
\end{equation*}
\end{Prop}

\noindent\underline{Proof}: The condition \eqref{C1} is equivalent to
$$\phi(vm_i,m_j)=\phi(m_i,vm_j),\qquad \forall\;0\leq i,j\leq N.$$
Let $(i,j)\in \{0,\dots,N\}^2$. Since $m_j\in V$, he have 
$$\phi(vm_i,m_j)=\phi(M(vm_i),m_j).$$
Moreover, \eqref{cond_moments_m_i} implies 
\begin{equation}\label{formule M(vm_i)}
M(vm_i)=m_{i-1}+a_i m_N,\qquad 0\leq i\leq N, 
\end{equation}
with the convention $m_{-1}=0$.
Thus
\[\phi(vm_i,m_j)=\phi(m_{i-1},m_j)+a_i\phi(m_N,m_j)=Q_{i-1,j}+a_i Q_{N,j}.
\]
But
$$(A^{T}Q)_{i,j}=\sum_{k=0}^N A_{k,i}Q_{k,j}=Q_{i-1,j}+a_i Q_{N,j},$$
with the convention $Q_{-1,j}=0$.
Therefore, we have $$\phi(vm_i,m_j)=(A^{T}Q)_{i,j}.$$
Thus, \eqref{C1} amounts to the matrix $A^{T}Q$ be (real and) symmetric, i.e.
$A^{T}Q=QA,$
since $Q$ is symmetric.
\cqfd
This result suggests a particular way of defining the inner product on $V$:
\begin{Cor}\label{choix_Q}
Choose $\phi$ s.t.
\begin{equation*}%\label{choix_Q}
Q=(P^{-1})^T\,P^{-1}=(P P^{T})^{-1},
\end{equation*}
where $P$ is defined by
$\eqref{matrice_de_passage}$.  
This choice implies 
\begin{equation}\label{valeurs du ps sur V}
\phi(m_i,m_j)=Q_{i,j}=\sum_{k=0}^N P_{k,i}^{-1}P_{k,j}^{-1}.
\end{equation}
\end{Cor}

\noindent \underline{Proof.}
We have, denoting by $D=diag(\lambda_0,\dots,\lambda_N)$:
\[\begin{split}
A^{T}Q=&A^{T}(P^{-1})^{T}P^{-1}=(P^{-1}A)^{T}P^{-1}
=(DP^{-1})^{T}P^{-1}=(P^{-1})^{T}(DP^{-1})\\
&=(P^{-1})^{T}P^{-1}A=QA,
\end{split}\]
thus $Q$ satisfies $A^TQ=QA$. Moreover, it is easy to check
that $Q$ is symmetric, definite, and positive.
\cqfd

\bigskip

In the rest of the proof, we always choose $\Phi$ according to \eqref{valeurs du ps sur V}. 
%%%%%%%%%%%%%%%%%%%%%%%%%%%%%%%%%%%%%%%%%%%%%%%%%%%%%%%%%%%%%%%%%%%%%%%%%%%%%%%%%%%%%%%%%%%%%%%%%%%%%%%%%%%%%%%%%%%%%%%%%%%%%%%%%%%%%%%%%%%%%%
\subsection{Study of the condition \eqref{C2}}
%%%%%%%%%%%%%%%%%%%%%%%%%%%%%%%%%%%%%%%%%%%%%%%%%%%%%%%%%%%%
Now, we assume \eqref{valeurs du ps sur V} to be satisfied 
and analyze \eqref{C2}.
\begin{Prop} Assume \eqref{valeurs du ps sur V}, and consider the
  following polynomials 
\begin{equation}\label{formule explicite des T_i}
T_i(X)=\frac{1}{Q_{0,N}}\sum_{k=0}^NQ_{k,i}X^k,\qquad 0\leq i\leq N, 
\end{equation}
where the matrix $Q\in M_{N+1}(\r)$ is defined by \eqref{choix_Q}.

Then, the condition \eqref{C2} is equivalent to
\begin{equation}\label{C2 condition 1}
 T_N(v)m_i(v)=T_i(v)m_N(v),\qquad v\in I,\qquad 0\leq i\leq N-1,
\end{equation}
\begin{equation}\label{C2 condition 2}
\forall f\in K,\qquad \phi_{|K}\left(\frac{1}{Q_{0,N}}
\times\frac{\chi_A(v)m_N(v)}{T_N(v)},f\right)=\int_{I}v^{N+1}f(v)dv.
\end{equation}\label{propc2}
\end{Prop}

The proof of this proposition is split in several lemmas. First we find two equivalent conditions to \eqref{C2}.
\begin{Lemma} The condition \eqref{C2} is equivalent to
\begin{equation}
\label{C2.1}
\forall i\in\{1,\dots,N\},\quad
vm_i-m_{i-1}-a_im_N=\frac{Q_{i,N}}{Q_{0,N}}(vm_0-a_0m_N)
\end{equation}
\begin{equation}
\label{C2.2}
\forall f\in K,\quad
\phi_{|K}\left(\frac{1}{Q_{0,N}}(vm_0-a_0m_N),f\right)=\int_{I}w^{N+1}f(w)dw,
\end{equation}\label{lem1}
\end{Lemma}
We then have to study conditions \eqref{C2.1} and \eqref{C2.2}. 

\noindent\underline{Proof of Lemma \ref{lem1}}: The condition \eqref{C2} is obviously equivalent to
$$\phi(vm_i,f)=\phi(m_i,vf),\qquad 0\leq i\leq N,\quad f\in K.$$
Let $i\in \{0,\dots,N\}$ and $f\in K$. Since $f\in K$, he have 
$$\phi(vm_i,f)=\phi_{|K}(vm_i-M(vm_i),f)=\phi_{|K}(vm_i-m_{i-1}-a_im_N,f),$$
using \eqref{formule M(vm_i)}.
Similarly, since $m_i\in V$, 
$$\phi(m_i,vf)=\phi_{|V}(m_i,M(vf))=Q_{i,N}\int_{I}w^{N+1}f(w)dw,$$
where the last equality is deduced from
\begin{equation*}%\label{formule M(vf)}
\forall f\in K,\qquad M(vf)=\left(\int_{I}w^{N+1}f(w)dw\right)m_N, 
\end{equation*}
recalling that $\int w^n f\,dv=0$, $\forall n\leq N$.

Thus, the condition \eqref{C2} amounts to for any $0\leq i\leq N$ and $f\in K$ 
\begin{equation*}%\label{C2 equiv}
\phi_{|K}(vm_i-m_{i-1}-a_im_N,f)=Q_{i,N}\int_{I}w^{N+1}f(w)dw.
\end{equation*}
from which we deduce \eqref{C2.2}, and
\[
\phi_{|K}\left(\frac{1}{Q_{i,N}}(vm_i-m_{i-1}-a_im_N)
-\frac{1}{Q_{0,N}}(vm_0-a_0m_N),f\right)=0.
\]
Therefore, for all $i\in \{1,\dots,N\}$, we have 
$$\frac{1}{Q_{i,N}}(vm_i-m_{i-1}-a_im_N)
-\frac{1}{Q_{0,N}}(vm_0-a_0m_N)\in K\cap K^{\bot,\phi}=\{0\},$$ 
which shows the relation \eqref{C2.1}.\\
Conversely, the conditions \eqref{C2.1}+\eqref{C2.2} imply
\eqref{C2} just by following the previous steps in reverse order. 
\cqfd

We start with condition \eqref{C2.1}
\begin{Lemma} Consider the polynomial
\begin{equation*}%\label{def D}
D(X)=\sum_{i=0}^N \beta_i X^i, 
\end{equation*}
where
\begin{equation*}%\label{def beta_i}
\beta_i=\frac{Q_{i,N}}{Q_{0,N}}=\frac{\phi(m_i,m_N)}{\phi(m_0,m_N)},\qquad
0\leq i\leq N.  
\end{equation*}
The condition \eqref{C2.1} is equivalent to 
\begin{equation*}%\label{C21 equiv}
D(v)m_i(v)=T_i(v)m_N(v),\qquad v\in I,\qquad 0\leq i\leq N,
\end{equation*}
where $(T_i)_{0\leq i\leq N}$ are the following polynomials:
\begin{equation}\label{formule de recurrence des T_i}
\begin{array}[]{llll}
T_N(X)=D(X),\\\\
T_{i-1}(X)=XT_i(X)-a_iD(X)-\beta_i \chi_A(X),\qquad 1\leq i\leq N.
\end{array}
\end{equation}\label{lemmac2.1}
\end{Lemma}

\noindent\underline{Proof.}
Setting $\beta_i=\frac{Q_{i,N}}{Q_{0,N}}$ for $i\in\{0,\dots,N\}$, we
deduce from \eqref{C2.1} the recursive formula 
\begin{equation*}%\label{formule recursive2 des m_i}
m_{i-1}=vm_i-(a_i-\beta_ia_0)m_N-\beta_i (vm_0),\qquad i=1\dots,N,
\end{equation*}
which leads to, for any $0\leq i\leq N$
\begin{equation*}%\label{formule des m_i en fonction de m_0 et m_N}
m_i=\left(v^{N-i}-\sum_{k=0}^{N-i-1}v^k(a_{i+k+1}-\beta_{i+k+1}
a_0)\right)m_N-\left(\sum_{k=1}^{N-i}\beta_{i+k}v^k\right)m_0.
\end{equation*}
Thus 
$$D(v)m_0(v)=\left(\ds v^N-\sum_{k=0}^{N-1}(a_{k+1}-\beta_{k+1}
a_0)v^k \right) m_N(v),$$
where
$$D(v):=\beta_0+\beta_1 v+\dots +\beta_N v^N.$$
We deduce 
$$vD(v)m_0(v)=\left(\ds \chi_A(v)+a_0D(v) \right) m_N(v).$$
Thus, it comes from the recursive formula on the $m_i$ that
\begin{equation*}%\label{formule recursive3 des m_i}
D(v)m_{i-1}(v)=vD(v)m_i(v)-(a_i D(v)+\beta_i \chi_A(v))
m_N(v),\qquad i=1\dots,N,
\end{equation*}
which is exactly
\begin{equation}\label{formule des m_i}
D(v)m_i(v)=T_i(v)m_N(v),\qquad i=0,\dots,N,
\end{equation}
with
\[
\begin{array}[]{llll}
T_N(X)=D(X),\\\\
T_{i-1}(X)=XT_i(X)-a_iD(X)-\beta_i \chi_A(X),\quad i=1,\dots N.
\end{array}\]
Conversely, if the functions $(m_i)_{0\leq i\leq N}$ satisfy these
relations, then \eqref{C2.1} is obvious, almost  
everywhere $v\in I$ (except at the roots of the polynomial $D$).
\cqfd

\begin{Rk} We notice that
\begin{equation}\label{lien entre psi_0 et P_A}
XT_0(X)-a_0D(X)=\chi_A(X). 
\end{equation}

Furthermore, the formula \eqref{formule de recurrence des T_i} implies:
$$1+deg(T_{i})\leq max(deg(T_{i-1}),N+1),\qquad 1\leq i\leq N.$$

But $T_0\in \r_N[X]$, thus $T_i\in\r_N[X]$ for all $i\in \{0,\dots,N\}$.
\end{Rk}
We can now prove
\begin{Lemma} We have the explicit formula \eqref{formule explicite des T_i} for the 
  $(T_i)_{0\leq i\leq N}$ which is recalled here
\[
T_i(X)=\frac{1}{Q_{0,N}}\sum_{k=0}^NQ_{k,i}X^k,\qquad 0\leq i\leq N. 
\]\label{lem2}
\end{Lemma}

\noindent\underline{Proof}: First, according to the definition of $D$ in Lemma \ref{lemmac2.1}, we have, for all
$k\in \{0,\dots,N\}$, 
\[\begin{split}
D(\lambda_k)=&\frac{1}{Q_{0,N}}\sum_{i=0}^N 
Q_{i,N}\lambda_{k}^i=\frac{1}{Q_{0,N}}\sum_{i=0}^N\sum_{j=0}^N 
P_{j,i}^{-1}P_{j,N}^{-1}P_{i,k}
=\frac{1}{Q_{0,N}}P_{k,N}^{-1}\\
=&\frac{1}{Q_{0,N}\,\,\pi_k}.
\end{split}\]
Then, the recursive formula \eqref{formule de recurrence des T_i}
easily implies: 
$$T_j(\lambda_k)=(\lambda_k^{N-j}-a_N\lambda_k^{N-j-1}
-\dots-a_{j+2}\lambda_k-a_{j+1})D(\lambda_k),\qquad 0\leq j\leq N$$
which gives (according to
\eqref{matrice_de_passage_inverse}-\eqref{autre formule inverse P}) 
\begin{equation}\label{valeurs des Ti}
 T_j(\lambda_k)=\tilde{p}_{k,j}D(\lambda_k)=
\frac{\tilde{p}_{k,j}}{Q_{0,N}\,\,\pi_k}=\frac{P_{k,j}^{-1}}{Q_{0,N}}
,\qquad 0\leq j,k\leq N.
\end{equation}
Since each polynomial $T_i$ have a degree equal to $N$, it is enough
to check equality \eqref{formule explicite des T_i}  on the set
$\{\lambda_0,\dots, \lambda_N\}$,  which allows to conclude.
\cqfd

\begin{Rk} The formula \eqref{formule explicite des T_i} shows that
  the polynomials $(T_i)_{0\leq i\leq N}$ 
 are a basis of $\r_N[X]$ since the matrix $Q$ is invertible.  
\end{Rk}
We may finally characterize \eqref{C2.2}
\begin{Lemma} We assume that \eqref{C2.1} holds. Then, the condition
  \eqref{C2.2} is equivalent to  
\begin{equation*}%\label{C22 equiv}
\forall f\in K,\qquad
\phi_{|K}\left(\frac{1}{Q_{0,N}}\times\frac{\chi_A(v)m_N(v)}
{T_N(v)},f\right)=\int_{I}v^{N+1}f(v)dv. 
\end{equation*}\label{lem3}
\end{Lemma}

\noindent\underline{Proof}: It is straightforward using \eqref{formule des m_i} and the
formula \eqref{lien entre psi_0 et P_A}. 
\cqfd

We now have all what is needed to prove Prop. \ref{propc2}

\noindent\underline{Proof of Prop. \ref{propc2}}:
By Lemma \ref{lem1}, condition \eqref{C2} is equivalent to \eqref{C2.1}-\eqref{C2.2}. By Lemma \ref{lemmac2.1}, condition \eqref{C2.1} is equivalent to \eqref{C2 condition 1} provided that the $T_i$ are defined by \eqref{formule de recurrence des T_i}. Lemma \ref{lem2} shows that the recursive formula \eqref{formule de recurrence des T_i} actually gives the explicit formula \eqref{formule explicite des T_i}. Finally by Lemma \ref{lem3}, we know that condition \eqref{C2.2} is equivalent to \eqref{C2 condition 2} thus concluding the proof.
\cqfd 
%%%%%%%%%%%%%%%%%%%%%%%%%%%%%%%%%%%%%%%%%%%%%%%%%%%%%%%%%%%%%%%%%%%%%%%%%%%%%%%%%%%%%%%%%%%%%%%%%%%%%%%%%%%%%%
%%%%%%%%%%%%%%%%%%%%%%%%%%%%%%%%%%%%%%%%%%%%%%%%%%%%%%%%%%%%%%%%%%%%%%%%%%%%%%%%%%%%%%%%%%%%%%%%%%%%%%%%%%%%%%%%%%%%%%%%%%%%%%%%
\subsection{Study of the condition \eqref{C3}}
%%%%%%%%%%%%%%%%%%%%%%%%%%%%%%%%%%%%%%%%%%%%%%%%%%
We prove
\begin{Prop} Assume \eqref{valeurs du ps sur V}+\eqref{C2 condition
    1}+\eqref{C2 condition 2}. Then, setting 
\begin{equation}\label{def_poids}
\rho(v):=\frac{Q_{0,N}T_N(v)}{m_N(v)},
\end{equation}
and
\begin{equation*}
\phi_{|K}(f,g)=\int_{I}f(v)g(v)\rho(v)dv,\qquad \forall (f,g)\in K^2, 
\end{equation*}
the condition \eqref{C3} is satisfied.\label{propphiK}
\end{Prop}

\noindent
\underline{Proof}: This choice is compatible with \eqref{C2.2}, since
\eqref{C2.2} can also be written as
$$\forall f\in K,\qquad \phi_{|K}\left(\frac{\chi_A}{\rho},f\right)
=\int_{I}v^{N+1}f(v)dv,$$
and since we have $\int_{I}\chi_A(v)f(v)dv=\int_{I}v^{N+1}f(v)dv$
since $f\in K$. 

Moreover, \eqref{C3} is obviously satisfied with this choice.
\cqfd
%%%%%%%%%%%%%%%%%%%%%%%%%%%%%%%%%%%%%%%%%%%%%%%%%%%%%%%%%%%%%
%%%%%%%%%%%%%%%%%%%%%%%%%%%%%%%%%%%%%%%%%%%%%%%%%%%%%%%%%%%%%%%
\subsection{Choice of $\phi$ on the subspace $V$} 
%%%%%%%%%%%%%%%%%%%%%%%%%%%%%%%%%%%%%%%%%%%%%%%%%%%%%%%%%%%%%%%%%
%%%%%%%%%%%%%%%%%%%%%%%%%%%%%%%%%%%%%%%%%%%%%%%%%%%%%%%%%%%%%%%%%
For the moment $\Phi$ is defined as a weighted $L^2$ type inner product on $K$ by Prop. \ref{propphiK} and on $V$ by Corollary \ref{choix_Q}.

We wish to define $\Phi$ as a weighted inner product on the whole $K\oplus V$. 
The following lemma shows that provided $\rho$ satisfies the right relations then $\Phi$ as defined by Prop. \ref{propphiK} and Corollary \ref{choix_Q} is automatically of the right form.
\begin{Lemma} Assume \eqref{valeurs du ps sur V}, \eqref{C2
    condition 1}, \eqref{C2 condition 2}, \eqref{def_poids},  
and assume that the weight $\rho$ satisfies the moment conditions:
\begin{equation*}%\label{premiere conditions de moments sur le poids}
\int_{I}\frac{R(v)}{\rho(v)}dv=\sum_{k=0}^N R(\lambda_k),\qquad R\in\r_{2N}[X].
\end{equation*}
Then, we have:
\begin{equation}\label{valeurs de phi sur V}
 \phi_{|V}(f,g)=\int_{I}f(v)g(v)\rho(v)dv.
\end{equation}\label{phiV}
\end{Lemma}
 
\noindent\underline{Proof}: It is sufficient to show that the formula
\eqref{valeurs de phi sur V} holds for $(f,g)=(m_i,m_j)$.  
We have, for $(i,j)\in \{0,\dots,N\}^2$,
$$\int_{I}m_i(v)m_j(v)\rho(v)dv=\int_{I}
\left(\tilde{T_i}\tilde{T_j}\right)(v)\times\frac{dv}{\rho(v)},$$ 
setting 

\begin{equation*}%\label{normalisation des Ti}
 \tilde{T_i}=Q_{0,N}T_i,\qquad 0\leq i\leq N.
\end{equation*}
Moreover, we have
$$\phi(m_i,m_j)=Q_{i,j}=\sum_{k=0}^N P_{k,i}^{-1}P_{k,j}^{-1}=\sum_{k=0}^N 
\left(\tilde{T_i}\tilde{T_j}\right)(\lambda_k),$$
according to \eqref{valeurs du ps sur V} and \eqref{valeurs des Ti}.

Since the $\frac{(N+1)(N+2)}{2}$ polynomials
$(\tilde{T_i}\tilde{T_j})_{0\leq i\leq j\leq N}$ are in $\r_{2N}[v]$,
we see 
that the assumption on $\rho$
guarantees that \eqref{valeurs de phi sur V} is satisfied.
\cqfd
%
%%%%%%%%%%%%%%%%%%%%%%%%%%%%%%%%%%%%%%%%%%%%%%%%%%%%%%%%%%%%%%%%%%%%%%%%%%%%%%%%%%%%%%%%%%%%%%%%%%%
%%%%%%%%%%%%%%%%%%%%%%%%%%%%%%%%%%%%%%%%%%%%%%%%%%%%%%%%%%%%%%%%%%%%%%%%%%%%%%%%%%%%%%%%%%%%%%%%%%%
\subsection{Proof of the theorem
\eqref{th maxwellienne}: Synthesis}
%%%%%%%%%%%%%%%%%%%%%%%%%%%%%%%%%%%%%%%%%%%%
We summarize here all the definitions and check rigorously that they are compatible.

So, assume there exists a function $\rho=\rho(v)$, positive on $I$ and satisfying
the moment conditions: 
 \begin{equation}\label{cond_moments_synthese}
\int_{I}\frac{R(v)}{\rho(v)}dv=\sum_{k=0}^{N}R(\lambda_{k}),\qquad
\forall R\in\r_{2N+1}[X].	 
\end{equation}
Note that this in particular implies that $1/\rho$ is integrable on $I$.

The space $E=L^2(I,\rho(v)dv)$ is a Hilbert space, for the inner product 
\[
\phi(f,g)=\int_{I}f(v)g(v)\rho(v)dv,\qquad (f,g)\in E^2.
\]
We define the map $M:E\rightarrow E$ by
\[
Mf=\sum_{i=0}^N\left(\int_{I}v^{i}f(v)dv\right)
\frac{\tilde{T_i}(v)}{\rho(v)},\qquad
f\in E,
\] 
where
$$\tilde{T_i}(X)=\sum_{k=0}^NQ_{k,i}X^k,\qquad 0\leq i\leq N,$$
and $Q=(P^{T})^{-1}P^{-1}$ (the matrix $P$ is defined by
\eqref{matrice_de_passage}). 
\begin{itemize}
\item[$\bullet$] First, the map $M$ is well defined as 
$\forall f\in E,\; \forall i\in\{0,\dots,N\}$
\[ \left|\int_{I}v^{i}f(v)dv\right|{}
\leq \left(\int_{I}\frac{|v|^{2i}}{\rho(v)}dv\right)^{1/2}
\left(\int_{I}|f(v)|^{2}\rho(v)dv\right)^{1/2},
\]
and thus $Mf$ makes sense. Moreover, for $f\in E$, we have $Mf\in E$ because 
\[
\begin{split}
&\int_{I}|Mf(v)|^{2}\rho(v)dv  = 
\int_{I}\left|\sum_{i=0}^{N}\left(\int_{I}w^{i}f(w)dw\right)\frac{\tilde
  T_{i}(v)}
{\rho(v)}\right|^{2}\rho(v)dv\\
 &\quad  \leq    \ds\int_{I}\left(\sum_{i=0}^{N}
\left|\int_{I}w^{i}f(w)dw\right|^{2}\right)
\left(\sum_{i=0}^{N}\left|\frac{\tilde T_{i}(v)}
{\rho(v)}\right|^{2}\right)\rho(v)dv\\
 &\quad \leq 
\ds\left(\sum_{i=0}^{N}\left|\int_{I}w^{i}f(w)dw\right|^{2}\right)
\times  
 \sum_{i=0}^{N}\int_{I}\frac{|\tilde T_{i}(v)|^{2}}{\rho(v)}dv\quad <\,\,\infty.
\end{split}
\]
\item[$\bullet$] It is obvious that $M$ is linear. We check that $M$
  is a projection: let $f\in E$, we have, for $0\leq j \leq N$, 
$$\begin{array}{lll}
\ds\int_{I}v^{j}Mf(v)dv & = & \ds\sum_{i=0}^{N}
\left(\int_{I}w^{i}f(w)dw\right)\left(\int_{I}\frac{v^{j}
\tilde T_{i}(v)}{\rho(v)}dv\right)\\\\
& = &  \ds\sum_{i=0}^{N}\left(\int_{I}w^{i}f(w)dw\right)
\left(\sum_{k=0}^{N}\lambda_{k}^{j}\tilde T_{i}(\lambda_{k})\right).
\end{array}
$$
Moreover, we have 
$$\lambda_{k}^{j}=P_{j,k},\qquad \tilde
T_{i}(\lambda_{k})=P_{k,i}^{-1},\qquad 0\leq i,j,k\leq N,$$ 
thus $$\sum_{k=0}^{N}\lambda_{k}^{j}\tilde
T_{i}(\lambda_{k})=\delta_{i,j},\qquad 0\leq i,j\leq N.$$ 
We deduce
$$\begin{array}{lll}
\ds\int_{I}v^{j}Mf(v)dv & = & \ds\int_{I}w^{j}f(w)dw,\qquad 0\leq j\leq N,
\end{array}
$$
which implies $M\circ M=M$.
\item[$\bullet$] $M$ is an orthogonal projection (for the inner
  product $\phi$), because it is a self-adjoint projector: 
\[\begin{split}
\phi(Mf,g) & =  \int_{I}Mf(v)g(v)\rho(v)dv\\
&=  
\sum_{i=0}^{N}\left(\int_{I}w^{i}f(w)dw\right)
\left(\int_{I}\tilde T_{i}(v)g(v)dv\right)\\
 & = 
\ds\sum_{i=0}^{N}\sum_{k=0}^{N}Q_{k,i}
\left(\int_{I}w^{i}f(w)dw\right)
\left(\int_{I}v^{k}g(v)dv\right)\,\,=\,\,\phi(f,Mg), 
\end{split}
\]	
as the matrix $Q$ is symmetric.
\item[$\bullet$] The Maxwellian $M$ satisfies moment conditions:
\begin{equation*}
\left\{\begin{array}{llll}
\ds\int_{I}v^i Mf(v) dv & = & \ds\int_{I}v^i f(v) dv &  i=0,\dots,N,\\\\
\ds\int_{I}v^{N+1} Mf(v) dv & = &\ds \sum_{i=0}^N a_i \int_{I}v^i f(v) dv & 
\end{array}
\right..
\end{equation*}
In fact, the first conditions have been already established, and the
second ones result from the following formula, using
\eqref{cond_moments_synthese} and since $\chi_{A}(v)\tilde T_{i}(v)$
is a polynomial of degree $2N+1$, 
$$ \begin{array}{llll}
\ds\int_{I}\chi_{A}(v)Mf(v) & = & \ds\sum_{i=0}^{N}
\left(\int_{I}w^{i}f(w)dw\right)\int_{I}
\frac{\chi_{A}(v)\tilde T_{i}(v)}{\rho(v)}dv\\\\
 & =
&\ds\sum_{i=0}^{N}\left(\int_{I}w^{i}f(w)dw\right)
\sum_{k=0}^{N}\chi_{A}(\lambda_{k})\tilde 
T_{i}(\lambda_{k})=0. 
\end{array}$$
\item[$\bullet$] By Corollary \ref{choix_Q}, Prop. \ref{propphiK} and Lemma \ref{phiV}, the inner product $\Phi$ satisfies the assumptions of Lemma \ref{estimation_formelle}. This concludes the proof of Theorem \ref{th maxwellienne}.\cqfd
\end{itemize}	
%%%%%%%%%%%%%%%%%%%%%%%%%%%%%%%%%%%%%%%%%%%%%%%%%%%%%%%%%%%%%%%%%%%%%%%%%%%%%%%%%%%%%%%%%%%%%%%%%ù%
%%%%%%%%%%%%%%%%%%%%%%%%%%%%%%%%%%%%%%%%%%%%%%%%%%%%%%%%%%%%%%%%%%%%%%%%%%%%%%%%%%%%%%%%%%%%%%%%%%%%
\subsection{From Th. \ref{th maxwellienne} to Th. \ref{th stabilite}: Uniform stability estimate on the BGK model}
%%%%%%%%%%%%%%%%%%%%%%%%%%%%%
%%%%%%%%%%%%%%%%%%%%%%%%%%%%%%%%%
The first point is to control the collision term which is done by	
\begin{Lemma} Define $\rho$ and $\Phi$ as per Th. \ref{th maxwellienne} then 
\begin{equation}\label{inegalite_lemme}
\begin{split}
\ds\int_{\r}\int_{I}f_{N}^{\eps}L(f_{N}^{\eps})\rho(v)dvdx  \leq &  \left(\Lambda_{N,d}	
\left(\int_{I}|q(v)|^{2}\rho(v)dv\right)^{1/2}-\lambda\right)\\
&\ \int_{\r}\int_{I}|f_{N}^{\eps}|^{2}\rho(v)dvdx,
\end{split}
\end{equation}
with 
\begin{equation*}%\label{def constante Lambda}
\Lambda_{N,d}=\sqrt{\sum_{j=0}^{d}\alpha_{j}^{2}}\sqrt{\sum_{j=0}^{d}\sum_{k=0}^{N}\lambda_{k}^{2j}}
=\sqrt{\sum_{j=0}^{d}\alpha_{j}^{2}}\sqrt{\sum_{k=0}^{N}\frac{1-\lambda_{k}^{2d+2}}{1-\lambda_{k}^{2}}}.	
\end{equation*}\label{controlL}
\end{Lemma}

\noindent\underline{Proof}: We have
\[\begin{split}
&\int_{\r}\int_{I}f_{N}^{\eps}(t,x,v)L(f_{N}^{\eps})(t,x,v)\rho(v)dvdx\\
&\ =  \int_{\r}\int_{I}f_{N}^{\eps}(t,x,v)q(v)\sum_{j=0}^{d}\alpha_{j}\left(\int_{I}w^{j}f_{N}^{\eps}(t,x,w)dw\right)\rho(v)dvdx\\
&\qquad
-\lambda\int_{\r}\int_{I}|f_{N}^{\eps}(t,x,v)|^{2}\rho(v)dvdx.\\
\end{split}
\]
So
\[\begin{split}
&\int_{\r}\int_{I}f_{N}^{\eps}(t,x,v)L(f_{N}^{\eps})(t,x,v)\rho(v)dvdx\\
&\ =  \int_{\r}dx\left(\ds\sum_{j=0}^{d}\alpha_{j}\int_{I}w^{j}f_{N}^{\eps}(t,x,w)dw\right)
\left(\int_{I}f_{N}^{\eps}(t,x,v)q(v)\rho(v)dv\right)\\
&\qquad
-\lambda\int_{\r}\int_{I}|f_{N}^{\eps}(t,x,v)|^{2}\rho(v)dvdx.
\end{split}\]
We can control the moments of $f_{N}^{\eps}$ in the following way:
\[
\left|\int_{I}w^{j}f_{N}^{\eps}(t,x,w)dw\right|^{2}\leq 
\left(\int_{I}\frac{w^{2j}}{\rho(w)}dw\right)\left(\int_{I}|f_{N}^{\eps}(t,x,w)|^{2}\rho(w)dw\right),
\]
and the assumption \eqref{conditions moments poids} implies 
\begin{equation*}%\label{controle des moments}
\left|\int_{I}w^{j}f_{N}^{\eps}(t,x,w)dw\right|^{2}
\leq \left(\sum_{k=0}^{N}\lambda_{k}^{2j}\right)\phi(f_{N}^{\eps},f_{N}^{\eps}).
\end{equation*}
We deduce that
\begin{equation*}%\label{controle_1}
\left|\sum_{j=0}^{d}\alpha_{j}\int_{I}w^{j}f_{N}^{\eps}(t,x,w)dw\right|^{2}\leq 
\left(\sum_{j=0}^{d}\alpha_{j}^{2}\right)
\left(\sum_{j=0}^{d}\sum_{k=0}^{N}\lambda_{k}^{2j}\right)\phi(f_{N}^{\eps},f_{N}^{\eps}).
\end{equation*}
Moreover, we have
\begin{equation*}%\label{controle_2}
\left|\int_{I}f_{N}^{\eps}(t,x,v)q(v)\rho(v)dv\right|^{2}=\phi(f_{N}^{\eps},q)^{2}\leq\phi_{N}(q,q)\,\phi(f_{N}^{\eps},f_{N}^{\eps}).	
\end{equation*}
Thus, combining the last two inequalities and using Cauchy-Schwarz inequality, we obtain
\[\begin{split}
\int_{\r}\!\!\int_{I}f_{N}^{\eps}L(f_{N}^{\eps})\rho(v)dvdx \leq&
 \sqrt{\sum_{j=0}^{d}\alpha_{j}^{2}}\sqrt{\sum_{j=0}^{d}\sum_{k=0}^{N}
\lambda_{k}^{2j}}\sqrt{\phi(q,q)}\,
\int_{\r}\phi(f_{N}^{\eps},f_{N}^{\eps})dx\\
&-\lambda\int_{\r}\phi(f_{N}^{\eps},f_{N}^{\eps})dx,
\end{split}\]
which is the desired estimate.
\cqfd

\medskip

Combining Th. \ref{th maxwellienne} and Lemma \ref{controlL}, we can obtain stability estimates for \eqref{BGK} uniform in $\eps$:
\begin{Prop} Let $f_{N}^{\eps}$ a solution of \eqref{BGK} with \eqref{initialdata}. Then, the following estimate holds for any $t>0$
\begin{equation}\label{estimation_de_stabilite_cinetique}
\sup_{\eps>0}\int_{\r}\int_{I}|f_{N}^{\eps}(t,x,v)|^{2}\rho(v)dvdx\leq 
e^{tC_{N,q}}\int_{\r}\int_{I}|f^{0}(x,v)|^{2}\rho(v)dvdx,
\end{equation}	
where 
\begin{equation*}%\label{constante_de_stabilite}
C_{N,q}=2\left(\ds\Lambda_{N,d}	
\left(\int_{I}|q(v)|^{2}\rho(v)dv\right)^{1/2}-\lambda\right),
\end{equation*}
and we recall
\[%\label{rappel constante Lambda}
\Lambda_{N,d}=\sqrt{\sum_{j=0}^{d}\alpha_{j}^{2}}\sqrt{\sum_{j=0}^{d}\sum_{k=0}^{N}\lambda_{k}^{2j}}
=\sqrt{\sum_{j=0}^{d}\alpha_{j}^{2}}\sqrt{\sum_{k=0}^{N}\frac{1-\lambda_{k}^{2d+2}}{1-\lambda_{k}^{2}}}.	
\]\label{BGKuniform}
\end{Prop}

\noindent \underline{Proof}:
It is straightforward using Gronwall lemma, as according to \eqref{bornes_stab_BGK} and \eqref{inegalite_lemme},
we have 
\[\begin{split}
\frac{d}{dt}\int_{\r}\int_{I}|f_{N}^{\eps}(t,x,v)|^{2}\rho(v)dvdx\leq 
&2\left(\ds\Lambda_{N,d}	
\left(\int_{I}|q(v)|^{2}\rho(v)dv\right)^{1/2}-\lambda\right)\\
&\int_{\r}\int_{I}|f_{N}^{\eps}(t,x,v)|^{2}\rho(v)dvdx.
\end{split}
\] 
\cqfd

\medskip

We are now ready to conclude the proof of Th. \ref{th stabilite}
Now, we show that we can pass to the limit $\eps\rightarrow 0$ in the BGK model \eqref{BGK}
to obtain a stability estimate on the hyperbolic system \eqref{systeme_hyperbolique}.

We fix $N\geq d$. Using \eqref{estimation_de_stabilite_cinetique}, we can see that the family
$(f_N^{\eps})_{\eps>0}$ is bounded in the space
$L^{2}(]0,+\infty[_{loc}\times \r_{x}\times I_{v},\rho(v)dtdxdv)$. Thus there exists a sequence 
$\eps_{k}\underset{k\rightarrow\infty}{\longrightarrow}0$ such that 
\[
f_{N}^{\eps_{k}}\underset{k\rightarrow\infty}{\tw} f_{N}
\in L^{2}\left(]0,+\infty[_{loc}\times \r_{x}\times I_{v},\rho_{N}(v)dtdxdv\right).
\]
Therefore 
\[
\del_{t}f_{N}^{\eps_{k}}+v\del_{x}f_{N}^{\eps_{k}}\underset{k\rightarrow\infty}{\tw} 
\del_{t}f_{N}+v\del_{x}f_{N},\qquad L(f_{N}^{\eps_{k}})\underset{k\rightarrow\infty}{\tw} L(f_{N})
\quad\mbox{in}\quad \mathcal{D}'(]0,\infty[\times \r\times I),
\]
from which we deduce 
\[
M(f_{N}^{\eps_{k}})-f_{N}^{\eps_{k}}=\eps_{k}(\del_{t}f_{N}^{\eps_{k}}+v\del_{x}f_{N}^{\eps_{k}}-L(f_{N}^{\eps_{k}}))
\underset{k\rightarrow\infty}{\tw} 0\quad\mbox{in}\quad \mathcal{D}'(]0,\infty[\times \r\times I),
\]
and thus
\[
M(f_{N}^{\eps_{k}})\underset{k\rightarrow\infty}{\tw} f_{N}\quad\mbox{in}\quad 
\mathcal{D}'(]0,\infty[\times \r\times I).
\]
Hence, passing to the limit in \eqref{conditions moments maxwellienne}, we show that the function $f_{N}$ is a 
kinetic interpretation of the system $\eqref{systeme_hyperbolique}$, in the sense that  
$\left(\mu_{i}=\int_{I}v^{i}f_{N}dv\right)_{0\leq i\leq N}$ satisfies \eqref{systeme_hyperbolique}. As this system is linear, hyperbolic, it has a unique solution for a given initial data. That means that any solution of \eqref{systeme_hyperbolique} can consequently be obtained as the moments of a limit $f_N$ of $f_N^\eps$. 

Moreover, $f_{N}$ also satisfies the bound \eqref{estimation_de_stabilite_cinetique} :
\begin{equation*}%\label{limite_estimation_de_stabilite_cinetique}
\forall t>0,\qquad \int_{\r}\int_{I}|f_{N}(t,x,v)|^{2}\rho(v)dvdx\leq 
e^{tC_{N,q}}\int_{\r}\int_{I}|f^{0}(x,v)|^{2}\rho(v)dvdx.
\end{equation*}
Thus we obtain the estimate \eqref{estimation_de_stabilite_macro}, since 
\[
\left|\mu_{i}(t,x)\right|^{2}=\left|\int_{I}v^{i}f_{N}(t,x,v)dv\right|^{2}
\leq \left(\sum_{k=0}^{N}\lambda_{k}^{2i}\right)\int_{I}|f_{N}(t,x,v)|^{2}\rho(v)dv,
\]
according to the assumptions on $\rho_{N}$. The proof of Th. \ref{th stabilite} is now complete.\cqfd
%%%%%%%%%%%%%%%%%%%%%%%%%%%%%%%%%%%%%%%%%%%%%%%%%%%%%%%%%%%%%%%%%%%%%%%%%%%%%%%%%%%%%%%%%%%%%%%%%%%%%%%%%%%%%%%%%%%%%%%%%%%%%%%%%%%%%%%%%%%%%%%%%%%%%%
%%%%%%%%%%%%%%%%%%%%%%%%%%%%%%%%%%%%%%%%%%%%%%%%%%%%%%%%%%%%%%%%%%
\section{Error estimate: Proof of Th. \ref{errorestimate} and Prop. \ref{resultat_smooth_case}}
%%%%%%%%%%%%%%%%%%%%%%%%%%%%%%%%%%%%%%%%%%%%%%%%%%%%%%%%%%%%%%%%%%%
%%%%%%%%%%%%%%%%%%%%%%%%%%%%%%%%%%%%%%%%%%%%%%%%%%%%%%%%%%%%%%%%%%%%%
%%%%%%%%%%%%%%%%%%%%%%%%%%%%%%%%%%%%%%%%%%%%%%%%%%%
\subsection{Estimates provided by the model}
%%%%%%%%%%%%%%%%%%%%%%%%%%%%%%%%%%%%%%%%%%%%%%%%%%%%%
In order to prove Th. \ref{errorestimate}, we will need good smoothness properties on the solution to the exact equation \eqref{eq_neutronics}. Fortunately this model is very simple to manipulate and the estimates we need easy to obtain in that case. 

Let us first start with the support in velocity
\begin{Prop} \label{support_vitesse_solution}
Assume \eqref{initialdata}, \eqref{forme_noyau} with $I\subset [-1,\ 1]$. Then the solution to \eqref{eq_neutronics} satisfies
\begin{equation*}\forall t\geq 0,\qquad a.e.\,\,x\in\r,\qquad {\rm supp_v}\, f(t,x,.)\subset I\subset [-1,\ 1].	
\end{equation*}	
\end{Prop}

\noindent\underline{Proof} :
If $f$ is a solution of \eqref{eq_neutronics} with $d=0$, then, using the Stokes formula, we have, formally   
\[
\begin{split}
\frac{d}{dt}\int_{\r}|f(t,x,v)|^{2}dx & =  2\int_{\r} f(t,x,v)q(v)\left(\int_{I}f(t,x,v^{*})dv^{*}\right)dx\\
&\qquad -2\lambda\int_{\r}|f(t,x,v)|^{2}dx\\	
& \leq   \ds 2q(v)\int_{\r} f(t,x,v)\left(\int_{I}f(t,x,v^{*})dv^{*}\right)dx.
\end{split}
\]
Therefore, integrating in time, we obtain
\[
\begin{split}
\int_{\r}|f(t,x,v)|^{2}dx  & \leq  \ds \int_{\r}|f^{0}(x,v)|^{2}dx \\
&\qquad + 2q(v)\int_{0}^{t}\int_{\r} f(s,x,v)\left(\int_{I}f(s,x,v^{*})dv^{*}\right)dxds,
\end{split}
\]
and since supp$\,q\subset I$, we get the result.
\cqfd

The model \eqref{eq_neutronics} also propagates the $H^k$-norm of the solution 
\begin{Prop}
Assume \eqref{initialdata}, \eqref{forme_noyau} with $I\subset [-1,\ 1]$. Then the solution to \eqref{eq_neutronics} satisfies
$$\forall t\geq 0,\qquad \|f(t)\|_{H^k(\r,\ L^{2}(I))}\leq e^{Ct}\|f^{0}\|_{H^k(\r,\ L^{2}(I))},$$
where $C=\|q\|_{L^{2}(I)}-\lambda.$	\label{hksolution}
\end{Prop}	

\noindent \underline{Proof}: First note that for any $k$, $\partial_x^k f$ is also a solution to \eqref{eq_neutronics} with $\partial_x^k f^0$ as initial data. Then 
\[\begin{split}
\frac{d}{dt}\|f(t)\|^{2}_{H^k(\r,\ L^{2}(I))}=&
2\sum_{p=0}^{k}\int_{\r}\left(\int_{\r}\del_{x}^{p}f(t,x,v)q(v)dv\right)
\left(\int_{I}\del_{x}^{p}f(t,x,v^{*})dv^{*}\right)\\
&-2\lambda \|f(t)\|^{2}_{H^k(\r,\ L^{2}(I))}.
\end{split}\]
Using Prop. \eqref{support_vitesse_solution} and H\"{o}lder inequality, we easily obtain 
\[
\frac{d}{dt}\|f(t)\|^{2}_{H^k(\r,\ L^{2}(I))}\leq 2\left(\|q\|_{L^{2}(I)}-\lambda\right)\|f(t)\|^{2}_{H^k(\r,\ L^{2}(I))},\]
and a simple Gronwall lemma gives the result.
\cqfd
 
We may finally conclude from Props. \ref{support_vitesse_solution}-\ref{hksolution} that
\begin{Cor}
Assume \eqref{initialdata}, \eqref{forme_noyau} with $I\subset [-1,\ 1]$. Then the moments of solution to \eqref{eq_neutronics} satisfy
\begin{equation}\label{propagation_moments}
\forall i\in\N,\quad\forall t\geq 0,\qquad \|\mu^{f}_{i}(t)\|_{L^{2}(\r)}\leq{}
e^{Ct}\|f^{0}\|_{L^{2}(\r^{2})},
\end{equation} 
and
\[
\sum_{i=0}^N \|\mu^{f}_{i}(t)\|_{L^{2}(\r)}\leq \, (2N+1)^{1/2}\,e^{Ct}\,\|f^0\|_{L^2(\r^2)}.
\]
where 
$C=\|q\|_{L^{2}(I)}-\lambda.$	\label{hkmoments}
\end{Cor} 

\noindent\underline{Proof.} Notice that
\[
|\mu_i^f|\leq \int_{-1}^1 |v|^i f(t,x,v)\,dv\leq \frac{1}{\sqrt{2i+1}}\,\left(\int_{-1}^1 |f(t,x,v)|^2\,dv\right)^{1/2}, 
\]
by Cauchy-Schwarz. Hence
\[
\|\mu^{f}_{i}(t)\|_{L^{2}(\r)}\leq \frac{1}{\sqrt{2i+1}}\,\|f(t)\|_{L^2(\r^2)},
\]
and one concludes by using Prop. \ref{hksolution}.\cqfd
%%%%%%%%%%%%%%%%%%%%%%%%%%%%%%%%%%%%%%%%%%%%%%%%%%%%%%%%%%%%%%%%%%
\subsection{Proof of Prop. \ref{resultat_smooth_case}: Error estimate in the smooth case}
%%%%%%%%%%%%%%%%%%%%%%%%%%%%%%%%%%%%%%%%%%%%%%%%%%%%%%%%%%%%%%%%%%%%	
%Here, we show how to control the error on the first moment, namely 
%\[
%\sup_{0\leq t\leq T}\|\mu_{0}^{f}(t,.)-\mu_{0}(t,.)\|_{L^{2}(\r)},
%\]
%in terms of the $N^{th}$ derivatives, in the case where the initial data is %smooth in the space variable $x$.

%\noindent\underline{Proof} : Since the initial data is smooth in $x$, it is easy to show that the solution $f=f(t,x,v)$ is 
%also $C^{\infty}$ in $x$.
%For $i\in\{0,\dots,N\}$, we set $\Delta_i=\mu_i-\mu_i^{f}$.
%We obviously have
%\begin{equation}\label{hierarchy_erreurs}
%\left\{\begin{aligned}
%&\del_t \Delta_i+\del_x \Delta_{i+1}  =  \ds\gamma_i
%\Delta_0-\lambda \Delta_i\qquad,  
%0\leq i\leq N,\\
%&\Delta_{N+1}(t,x)  =  \ds\sum_{i=0}^N a_i
%\mu_i- \mu_{N+1}^{f},
%\Delta_i(0,x)  =  0, \qquad 0\leq i\leq N.\\
%\end{aligned}\right.
%\end{equation}
%
%We can control the first error $\Delta_{0}$ in function of the last derivatives in 
%$x$, 
Since the functions $(\mu_{i})$ are smooth in the space variable, we have by \eqref{method_moments_d=0}, for all~$0\leq i\leq N$,
\[
\begin{split}
&\ds\frac{d}{dt}\|\mu_{i}(t)\|_{L^{2}(\r)}^{2}  =  2\ds\int_{\r}\mu_{i}(t,x)\del_{t}\mu_{i}(t,x)dx\\
&\quad =  -2\int_{\r}\mu_{i}(t,x)\del_{x}\mu_{i+1}(t,x)dx\\
&\quad + 2\gamma_{i}\int_{\r}\mu_{i}(t,x)\mu_{0}(t,x)dx  - 2\lambda\int_{\r}|\mu_{i}(t,x)|^{2}dx,\\ 
\end{split},
\]
or
\[
\begin{split}
&\ds\frac{d}{dt}\|\mu_{i}(t)\|_{L^{2}(\r)}^{2}
\leq  2\|\mu_{i}(t)\|_{L^{2}(\r)}\|\del_{x}\mu_{i+1}(t)\|_{L^{2}(\r)}\\
&\qquad\qquad+2|\gamma_{i}|\;\|\mu_{i}(t)\|_{L^{2}(\r)}\|\mu_{0}(t)\|_{L^{2}(\r)}
- 2\lambda\|\mu_{i}(t)\|_{L^{2}(\r)}^{2}.
\end{split}\]
Thus, we obtain
\begin{equation*}%\label{relation_entre_les_delta_i}
\begin{array}{lll}
\ds	\frac{d}{dt}\|\mu_{i}(t)\|_{L^{2}(\r)} & \leq & \|\del_{x}\mu_{i+1}(t)\|_{L^{2}(\r)}+|\gamma_{i}|\,\|\mu_{0}(t)\|_{L^{2}(\r)}
-\lambda\|\mu_{i}(t)\|_{L^{2}(\r)}.
\end{array}
\end{equation*}
Of course the same computation can be performed on the $\del_{x}^{k}\mu_{i}$, obtaining 
\begin{equation}\label{relation_entre_les_derivees_des_delta_i}
\begin{split}
&\forall k\in\N,\quad \forall i\in\{0,\dots,N\},\\
&\frac{d}{dt}\|\del_{x}^{k}\mu_{i}(t)\|_{L^{2}(\r)}  \leq  
\|\del_{x}^{k+1}\mu_{i+1}(t)\|_{L^{2}(\r)}+|\gamma_{i}|\;
\|\del_{x}^{k}\mu_{0}(t)\|_{L^{2}(\r)}\\
&\qquad\qquad\qquad\qquad-\lambda\|\del_{x}^{k}\mu_{i}(t)\|_{L^{2}(\r)}.
\end{split}
\end{equation}
This sequence of inequalities lets us control $\|\mu_{0}(t)\|_{L^{2}(\r)}$ in term of the "last 
derivatives", namely $\left(\|\del_{x}^{N}\mu_{i}\|_{L^{2}(\r)}\right)_{0\leq i\leq N}$. To do that, we set 
\[
\forall t\geq 0,\quad\forall k\in\{0,\dots,N\},\qquad H_{k}(t):=\sum_{i=0}^{k}\|\del_{x}^{k}\mu_{i}\|_{L^{2}(\r)}.
\]
The coefficients $\gamma_i$ are easily bounded by $\|q\|_{L^2}$. Hence the inequalities \eqref{relation_entre_les_derivees_des_delta_i}  imply  
\begin{equation}\label{relation_entre_H_k}
\forall t\geq 0,\quad\forall k\in\{0,\dots,N\},\quad
\ds	\frac{d}{dt}H_{k}(t) \leq  
H_{k+1}(t)+\left(2\|q\|_{L^{2}(I)}-\lambda\right)H_{k}(t).
\end{equation}
We now use  
\begin{Lemma}
Let $(H_{k}(t))_{k\geq 1}$ be a sequence of nonnegative $C^{1}$ functions such that:
$$\forall k\in\N,\quad \forall t\geq 0,\qquad \left\{\begin{array}{ll}
H_{k}'(t)\leq CH_{k}(t)+H_{k+1}(t)\\ H_{k}(0)=0	
\end{array}\right.,$$
where $C>0$ is a numerical constant independant of $k$.\\ 
Then, we have:
$$\forall p\in\N^{*},\quad \forall t\geq 0,\qquad 
H_{0}(t)\leq \frac{1}{(p-1)!}\int_{0}^{t}(t-s)^{p-1}e^{C(t-s)}H_{p}(s)ds.$$  	
\end{Lemma}

\noindent \underline{Proof}:
First, the assumption may be rewritten as
$$\forall k\in\N,\quad \forall t\geq 0,\qquad 
\frac{d}{dt}\left(H_{k}(t)e^{-Ct}\right)\leq H_{k+1}(t)e^{-Ct},$$
thus, integrating in $t$, we obtain:
$$\forall k\in\N,\quad \forall t\geq 0,\qquad 
H_{k}(t)\leq \int_{0}^{t}e^{C(t-s)}H_{k+1}(s)ds.$$ 
A simple recursion allows to conclude.
\cqfd

\bigskip

\noindent\underline{End of the proof of prop \eqref{resultat_smooth_case}}: Applying the previous lemma, we obtain 
$$\forall t\geq 0,\qquad H_{0}(t)\leq \frac{1}{(N-1)!}\int_{0}^{t}(t-s)^{N-1}e^{(2\|q\|_{L^{2}(I)}-\lambda)(t-s)}H_{N}(s)ds,$$ 
and thus, for all $t\in [0,T]$,
\[\begin{split}
&\|\mu_{0}\|_{L^{2}(\r)}  \leq  \ds\frac{1}{(N-1)!}\int_{0}^{t}(t-s)^{N-1}e^{(2\|q\|_{L^{2}(I)}-\lambda)(t-s)}\sum_{i=0}^{N}\|\del_{x}^{N}\mu_{i}(s)\|_{L^{2}(\r)}ds\\
&\quad \leq  \ds\frac{e^{(2\|q\|_{L^{2}(I)}-\lambda)T}}{(N-1)!}\int_{0}^{t}(t-s)^{N-1}\sum_{i=0}^{N}\|\del_{x}^{N}\mu_{i}(s)\|_{L^{2}(\r)}ds\\
%&\quad \leq \ds \frac{e^{(2\|q\|_{L^{2}(I)}-\lambda)T}}{(N-1)!}\sum_{i=0}^{N}\left(\int_{0}^{t}(t-s)^{2N-2}ds\right)^{1/2}
%\left(\int_{0}^{t}\|\del_{x}^{N}\Delta_{i}(s)\|_{L^{2}(\r)}^{2}\right)^{1/2}\\
%&\quad =  \ds\left(\frac{T^{2N-1}}{2N-1}\right)^{1/2}\frac{e^{(2\|q\|_{L^{2}(I)}-\lambda)T}}{(N-1)!}\sum_{i=0}^{N}\|\del_{x}^{N}\Delta_{i}\|_{L^{2}([0,T]\times \r)},
\end{split}\]
which gives
\begin{equation}\label{estimation_provisoire_smooth_case}
\sup_{0\leq t\leq T}\|\mu_{0}(t,.)\|_{L^{2}(\r)} \leq C \frac{T^{N}}{N!}\times
\sum_{i=0}^{N}\|\del_{x}^{N}\mu_{i}\|_{L^\infty_tL^{2}_x},
\end{equation}
where the constant $C$ depends on $T$, $\lambda$ and $q$. This concludes the proof of Prop. \ref{resultat_smooth_case}.\cqfd

%Furthermore, since the last derivative $\del_{x}^{N}f$ satisfy the same kinetic equation as $f$, we have, according to 
%\eqref{propagation_moments}, 
%$$\|\del_{x}^{N}\mu_{i}^{f}(t)\|_{L^{2}(\r)}=\|\mu_{i}^{\del_{x}^{N}f}(t)\|_{L^{2}(\r)}\leq 
%\left(\frac{1}{2}\right)^{i}e^{(\|q\|_{L^{2}(I)}-\lambda)t}\|\del_{x}^{N}f^{0}\|_{L^{2}(\r^{2})},$$ 
%thus 
%$$\|\del_{x}^{N}\mu_{i}^{f}\|_{L^{2}([0,T]\times\r)}\leq 
%\left(\frac{1}{2}\right)^{i}e^{(\|q\|_{L^{2}(I)}-\lambda)T}\sqrt{T}\|\del_{x}^{N}f^{0}\|_{L^{2}(\r^{2})},$$ 
%from which we deduce 
%\begin{equation}\label{estimation_provisoire_smooth_case_2}
%\sum_{i=0}^{N}\|\del_{x}^{N}\mu_{i}^{f}\|_{L^{2}([0,T]\times\r)}\leq{}
%2\sqrt{T}e^{(\|q\|_{L^{2}(I)}-\lambda)T}\|\del_{x}^{N}f^{0}\|_{L^{2}(\r^{2})}.
%\end{equation}
%
%Moreover, applying the stability estimate \eqref{estimation_de_stabilite_macro_d=0_sous_hypotheses} to $\left(\del_{x}^{N}\mu_{i}\right)_{0\leq i\leq N}$ 
%(which also satisfy the method \eqref{method_moments_d=0}),
%\begin{equation}\label{estimation_provisoire_smooth_case_3}
%\sum_{i=0}^{N}\|\del_{x}^{N}\mu_{i}\|_{L^{2}([0,T]\times\r)}\leq
%CN^{\gamma}\|\del_{x}^{N}f^{0}\|_{L^{2}(\r^{2})},
%\end{equation}
%where the constant $C\geq 0$ depends on $T$, $q$ and $\lambda$. 
%
%We conclude by gathering the estimates \eqref{estimation_provisoire_smooth_case}, \eqref{estimation_provisoire_smooth_case_2} and 
%\eqref{estimation_provisoire_smooth_case_3}.
%\cqfd
%
%
%%%%%%%%%%%%%%%%%%%%%%%%%%%%%%%%%%%%%%%%%%%%%%%%%%%
\subsection{Proof of Th. \ref{errorestimate}: Error estimate in the general case}
%%%%%%%%%%%%%%%%%%%%%%%%%%%%%%%%%%%%%%%%%%%%%%%%%%%
The general idea for the proof is to regularize the initial data. Then we have to bound 3 terms. First the error between the exact solution and the truncated hierarchy for this regularized initial data. This term is controlled by Prop. \ref{resultat_smooth_case}. The next term is the difference between the solution for the non regularized initial data and the solution for the regularized one, both solutions to the exact equation \eqref{eq_neutronics}. This is bounded using the estimates for \eqref{eq_neutronics}. The final term is the difference between the solution for the non regularized initial data and the solution for the regularized one, but both solutions to the truncated hierarchy \eqref{method_moments_d=0}. We bound this term thanks to assumption
\eqref{estimation_de_stabilite_macro_d=0_sous_hypotheses}.
 
\begin{itemize}
\item\textit{Step $1$: Regularization of the initial data} 

\noindent We fix $\eps>0$ and we choose $f_{\eps}^0\in
H^{k}(\r^2)$, with ${\rm
supp_v}\,f_{\eps}^0\subset I$, and such that (see the appendix for more details)
\begin{equation}\label{hypotheses sur la regularisee}
\begin{array}{lll}
a.e.\,\,v\in I, & f_{\eps}^0(.,v)\in C^{\infty}(\r),\\\\
& \|f_{\eps}^0(v)-f^0(v)\|_{L^2(\r)}\leq
C_k\eps^k\|f^0(v)\|_{H^k(\r)},\\\\
& \|f_{\eps}^0(v)\|_{H^m(\r)}\leq
C_k\eps^{k-m}\|f^0(v)\|_{H^k(\r)}\quad\forall m>k.
\end{array}
\end{equation}
Let $f_{\eps}=f_{\eps}(t,x,v)$ be the unique solution to Eq. \eqref{eq_neutronics} with $f_\eps(t=0)=f_\eps^0$.
%\begin{equation}\left\{\begin{array}{lll}
%\partial_t f_{\eps}+v\del_x f_{\eps}=L(f_{\eps})(t,x,v),\quad t\in \r_+,\
%(x,v)\in \r\times I,\\\\
% f_{\eps}(t=0,x,v)=f_{\eps}^0(x,v).
%\end{array}\right.
%\label{eqcin approchee}
%\end{equation}
%
We define the moments as usual
\begin{equation*}
%\label{def moments approches}
\mu_i^{f_{\eps}}(t,x)=\int_{I} v^i\,f_{\eps}(t,x,v)\,dv,\qquad i\in\N.
\end{equation*}
Of course, those moments still satisfy the hierarchy
\begin{equation*}
\partial_t \mu_i^{f_{\eps}}+\partial_x \mu_{i+1}^{f_{\eps}}=\ds\gamma_i
\mu_0^{f_{\eps}}-\lambda \mu_i^{f_{\eps}}\qquad
i\in \N.
%\label{eq moments approches}
\end{equation*}
We denote by $\mu_i^\eps$ the solution to the truncated hierarchy \eqref{method_moments_d=0} for the initial data 
%
%\begin{equation}%\label{method_moments_approchee_d=0}
%\left\{\begin{array}{lll}
%   \del_t \mu_i^{\eps}+\del_x \mu_{i+1}^{\eps}=\ds\gamma_i
%\mu_0^{\eps}-\lambda \mu_i^{\eps},\qquad i=0,\dots,N\\
%\mu_{N+1}^{\eps}=\ds\sum_{i=0}^N a_i \mu_i^{\eps}.
%  \end{array}\right.,
%\end{equation}
%
\begin{equation}
\mu_i^{\eps}(t=0,x)=\mu_i^{f^0_{\eps}}(x)=\int_I v^i f^0_{\eps}(x,v)\,dv,\qquad
0\leq i\leq N.
%\label{initialdata approchee}
\end{equation}

\smallskip

\item\textit{Step $2$: Error estimates in term of $\eps$.}
First note that $\mu_i^\eps-\mu_i^{f_\eps}$ satisfies the assumptions of Prop. \ref{resultat_smooth_case}. On the one hand
\[\begin{split}
\sum_{i\leq N}\|\partial_x^N (\mu_{i}^\eps-\mu_i^{f_\eps})\|_{L^2}&\leq \sum_{i\leq N} \left(\|\partial_x^N \mu_{i}^\eps\|_{L^2}+\|\partial_x^N \mu_{i}^{f_\eps}\|_{L^2}\right)\\
&\leq \sum_{i\leq N}\|\partial_x^N \mu_{i}^\eps\|_{L^2}+e^{Ct}\,(2N+1)^{1/2} \|\partial_x^N\,f^0_\eps\|_{L^2},
\end{split}\]
by Corollary \ref{hkmoments}. On the other hand by applying \eqref{estimation_de_stabilite_macro_d=0_sous_hypotheses} to $\mu_i^\eps$ and $0$, one gets
\[
\sum_{i\leq N}\|\partial_x^N\mu_{i}^\eps\|_{L^\infty_t L^2_x}\leq C\,N^\gamma\,\|\del_{x}^{N}f^0_\eps\|_{L^2_{x,v}}.
\]
So applying Prop. \ref{resultat_smooth_case}, we have for $\gamma'=\max(\gamma,1/2)$
\begin{equation}\label{formule erreur 1}
\begin{split}
\|\mu_0^{\eps}-\mu_0^{f_{\eps}}\|_{L^{\infty}([0,T],L^2(\r))}&\leq C\, 
\frac{N^{\gamma'}T^N}{N!}\times
\|\del_x^N f^{0}_{\eps}\|_{L^2(\r^{2})}\\
&\leq C_k\,\eps^{k-N}\, 
\frac{N^{\gamma'}T^N}{N!}\times
\|f^{0}\|_{H^k_x\,L^2_v},
\end{split}\end{equation}
by \eqref{hypotheses sur la regularisee}.

We can use the stability estimate \eqref{estimation_de_stabilite_macro_d=0_sous_hypotheses} to control 
%$\|\mu_0^{\eps}-\mu_0\|_{L^{\infty}([0,T],L^2(\r))}$
\begin{equation}\label{formule erreur 2}
\|\mu_0^{\eps}-\mu_0\|_{L^{\infty}([0,T],L^2(\r))}
%\sum_{i=0}^{N}\|\mu_{i}^{\eps}-\mu_{i}\|_{L^{\infty}([0,T],L^2(\r))}
\leq CN^{\gamma}\|f^{0}-f^{0}_{\eps}\|_{L^{2}}\leq C_k\,N^\gamma\,\eps^k\|f^{0}\|_{H^k_x\,L^2_v},
\end{equation}
again by \eqref{hypotheses sur la regularisee}.

At last, we can control $\|\mu_0^{f}-\mu_0^{f_{\eps}}\|_{L^{\infty}([0,T],L^2(\r))}$ according to Corollary \ref{hkmoments}
\begin{equation}\label{formule erreur 3}\begin{split}
\|\mu_0^{f}-\mu_0^{f_{\eps}}\|_{L^{\infty}([0,T],L^2(\r))}&=\|\mu_0^{f-f_{\eps}}\|_{L^{\infty}([0,T],L^2(\r))}\leq C\|f^{0}-f^{0}_{\eps}\|_{L^{2}(\r^{2})}\\
&\leq C_k\,\eps^{k}\,\|f^{0}\|_{H^k_x\,L^2_v}.
\end{split}\end{equation}
\item\textit{Step $3$: Choice of the parameter $\eps$}

We deduce from \eqref{formule erreur 1}, \eqref{formule erreur 2} and \eqref{formule erreur 3} the complete error estimate
\begin{equation}\label{erreur finale}
\|\mu_0-\mu_0^{f}\|_{L^{\infty}([0,T],L^2(\r))}\leq C\,
\left(N^{\gamma'}\,\frac{T^N}{N!}\,\eps^{k-N}+N^{\gamma}\,\eps^k\right)\,\|f^{0}\|_{H^k_x\,L^2_v},
\end{equation}
where the numerical constant $C\geq 0$ depends on $k$, $T$, $q$, and $\lambda$.
 
We of course choose the "best" value of $\eps$, 
which minimizes the error.
When $N>k$, we get for $\gamma\geq 1/2$
\[
\eps=\eps^{*}=\left(\frac{T^N}{N!}\right)^{1/N}\underset{N\rightarrow \infty}{\sim}
\frac{eT}{N}.
\]
If $\gamma<1/2$ then one takes instead
\[
\eps=\eps^{*}=\left(\frac{T^N\,N^{1/2-\gamma}}{N!}\right)^{1/N}\underset{N\rightarrow \infty}{\sim}
\frac{eT}{N},
\]
with the same asymptotic behaviour.

In both cases 
\begin{equation}\label{conclusion}
\begin{split}
 \|\mu_0-\mu_0^{f}\|_{L^{\infty}([0,T],L^2(\r))} & \leq 
\ds C_{T,q,\lambda,k}\,N^{\gamma}\,\left(\eps^{*}\right)^{k}\\
%\left(1+\frac{T^N N^{1/2}}{N!\left(\eps^{*}\right)^{N}}\right)\\ 
%& \leq  \ds C_{T,q,\lambda,k}\,N^{\gamma}\,\left(\eps^{*}\right)^{k}
%\left(\frac{N}{N-k}\right)\\\\
& \underset{N\rightarrow \infty}{\sim} 
\frac{C_{T,q,\lambda,k}}{N^{k-\gamma}},
\end{split}
\end{equation}
which concludes the proof of Theorem \ref{errorestimate}.\cqfd

\end{itemize}
%%%%%%%%%%%%%%%%%%%%%%%%%%%%%%%%%%%%%%%%%%%%
%%%%%%%%%%%%%%%%%%%%%%%%%%%%%%%%%%%%%%%%%%%
\section{Proof of Th. \ref{th stabilite tchebychev}: Example of the Tchebychev points\label{Tcheby}}
%%%%%%%%%%%%%%%%%%%%%%%%%%%%%%%%%%%%%%%%%%%
We now make the assumptions in Th. \ref{th stabilite tchebychev} and in particular that $d=0$ and that the $\lambda_k$ satisfy \eqref{tcheblambda}.

The constant $\Lambda_{N,0}$ defined by \eqref{constante Lambda} satisfies
\begin{equation}
\Lambda_{N,0}=(N+1)^{1/2}.\label{Lambda0}
\end{equation}
Moreover we define
\[
\rho_N(v)=\frac{\pi}{N+1}\,\sqrt{1-v^2}.
\]
In fact, the function $\frac{1}{\rho_{N}}$ is a normalization of the Tchebychev weight.	 

As the $\lambda_k$ gives the usual method of integration we have 
\begin{Prop}
We have, for all $N\geq 1$ and for all $R\in\r_{2N+1}[X]$ :
\begin{equation}\label{example_weight}
\int_{I}\frac{R(v)}{\rho_{N}(v)}dv=\frac{N+1}{\pi}\int_{]-1,1[}\frac{R(v)}{\sqrt{1-v^{2}}}dv
=\sum_{k=0}^N R(\lambda_k),\qquad \forall R\in\r_{2N+1}[X].
\end{equation}\label{rhotcheby}
\end{Prop} 
We may therefore apply Th. \ref{th stabilite} to this choice of $\lambda_k$ and for this choice of $\rho_N$.

Compute
\[
C_N(q)=\sqrt{N+1}\,\left(\int_{-1}^1 |q(v)|^2\rho_N(v)\,dv\right)^{1/2}\leq C\,\|q\|_{L^2}.
\]
Similarly
\[
C_N(f^0)\leq C\,N^{-1/2}\,\|f^0\|_{L^2}.
\]
So by Th. \ref{th stabilite}
\[
\sup _{i\leq N} \|\mu_i(t)\|_{L^2}\leq C\,e^{C\,T\,\|q\|_{L^2}}\,\|f^0\|_{L^2}\,N^{-1/2}\, \left(\sum_{k\leq N} \lambda_k^{2i}\right)^{1/2}\leq C\,e^{C\,T\,\|q\|_{L^2}}\,\|f^0\|_{L^2},
\]
which is exactly \eqref{estimation_de_stabilite_macro_tchebychev}. Moreover
\[
\sum_{i\leq N} \|\mu_i(t)\|_{L^2}\leq C\,\|f^0\|_{L^2}\,N^{-1/2}\,\sum_{i\leq N} \left(\sum_{k\leq N} \lambda_k^{2i}\right)^{1/2}.
\]
Of course by Cauchy-Schwarz, we have, for all $0\leq L\leq N$,
\[\begin{split}
N^{-1/2}\,\sum_{i\leq N} &\left(\sum_{k=L\ldots N-L} \lambda_k^{2i}\right)^{1/2}\leq \left(\sum_{i\leq N}\sum_{k=L\ldots N-L} \lambda_k^{2i}\right)^{1/2}\\
&\leq \left(\sum_{L\leq k\leq N-L} \frac{1}{1-|\lambda_k|^2}\right)^{1/2}
\end{split}\]
Note that
\[
1-|\lambda_k|^2\geq \frac{(k+1)^2}{C\,N^2}\quad \mbox{if}\ k\leq N/2,\quad 
1-|\lambda_k|^2\geq \frac{(N-k+1)^2}{C\,N^2}\quad \mbox{if}\ k\geq N/2.
\]
Hence, if $L\leq N/2$,
\[
N^{-1/2}\,\sum_{i\leq N} \left(\sum_{k=L\ldots N-L} \lambda_k^{2i}\right)^{1/2}\leq
C\,\left(\sum_{L\leq k\leq N/2} \frac{N^2}{(k+1)^2}\right)^{1/2}\leq C\,\frac{N}{L^{1/2}}. 
\]
Therefore
\[
\sum_{i\leq N} \|\mu_i(t)\|_{L^2}\leq C\,\|f^0\|_{L^2}\,(N^{1/2}\,L^{1/2}+N\,L^{-1/2}),
\]
and choosing $L=\sqrt{N}$ we obtain that this method satisfies the estimate \eqref{estimation_de_stabilite_macro_d=0_sous_hypotheses} with $\gamma=3/4$.
It only remains to apply Th. \ref{errorestimate} to conclude.\cqfd
%~ \begin{Rk}
%~ We can obtain a stability bound with a mild growth (polynomial) if the initial data $f^{0}$, the kernel $q$
%~ and the eigenvalues $(\lambda_{k})_{0\leq k\leq N}$ satisfy the 
%~ following assumptions 
%~ \begin{equation*}%\label{hypothese_donnee_initiale}
%~ \exists C\geq 0,\quad\exists\delta\geq 0,
%~ \qquad \int\!\!\!\int_{\r\times B(0,1/2)} |f^{0}(x,v)|^{2}\rho_{N}(v)dxdv\leq C N^{2\delta}\|f^{0}\|_{L^{2}(\r^{2})}^{2},
%~ \end{equation*}	
%~ %
%~ \begin{equation*}%\label{hypothese_noyau_q}
 %~ \exists C\geq 0,\quad\int_{B(0,1/2)}|q(v)|^{2}\rho_{N}(v)dv\leq \frac{C}{N},
%~ \end{equation*}
%~ %
%~ \begin{equation*}%\label{hypothese_valeurs_propres}
	 %~ \exists C\geq 0,\quad\exists\delta'\geq 0,\qquad\sum_{i=0}^{N}\left(\sum_{k=0}^{N}\lambda_{k}^{2i}\right)^{1/2}
	 %~ \leq C N^{\delta'}.
%~ \end{equation*}
%~ %
%~ \end{Rk}
%~ %
%%%%%%%%%%%%%%%%%%%%%%%%%%%%%%%%%%%%%%%%%%%%%%
%%%%%%%%%%%%%%%%%%%%%%%%%%%%%%%%%%%%%%%%%%
\section{Appendix}

The natural way to regularize is by convolution. However to obtain high order approximation, it is necessary to choose correctly the mollifier. 
In the $L^2$ framework though, things are quite simple by truncating in Fourier.
\begin{Prop} Let $k$ a positive integer, $f\in H^k(\r^d)$ and $\eps>0$. It
exists $f_{\eps}\in H^{\infty}(\r^d)$ such that
\begin{equation}\label{approx 1}
\|f-f_{\eps}\|_{L^2(\r^d)}\leq \eps^k \|D^k f\|_{L^2(\r^d)},
\end{equation}
\begin{equation}\label{approx 2}
\|D^m f_{\eps}\|_{L^2(\r^d)}\leq \eps^{k-m} \|D^k
f\|_{L^2(\r^d)},\quad\forall m\in\N.
\end{equation}

\end{Prop}

\noindent \underline{Proof} : We use Fourier's analysis. We consider $f_{\eps}\in
L^2(\r^d)$ defined by $$\widehat{f_\eps}(\xi)=\widehat{f}(\xi)\1_{\{|\xi|\leq
1/\eps\}}.$$
First, we have 
\[\begin{split}
\|\widehat{f}-\widehat{f_\eps}\|_{L^2(\r^d)}&=\left(\int_{
|\xi|>1/\eps } |\widehat{f}(\xi)|^2d\xi\right)^{1/2}\leq \left(\int_{
\r^d } \eps^{2k}
|\xi|^{2k}|\widehat{f}(\xi)|^2d\xi\right)^{1/2}\\
&=\eps^k\|\widehat{D^k
f}\|_{L^2(\r^d)}, 
\end{split}
\]
which proves \eqref{approx 1}.

On the other hand, for all $m\in\N$,
\[\begin{split}
\||\xi|^m\widehat{f_{\eps}}\|_{L^2(\r^d)}&=\left(\int_{
|\xi|\leq 1/\eps }|\xi|^{2m} |\widehat{f}(\xi)|^2d\xi\right)^{1/2}\!\!\leq
\left(\int_{\r^d}\eps^{2k-2m}|\xi|^{2k}|\widehat{f}(\xi)|^2d\xi\right)^{1/2}
\!\!\!\\
&=\eps^{k-m}\|\widehat{D^k f}\|_{L^2(\r^d)},
\end{split}\]
thus, $f_{\eps}\in H^m(\r^d)$ and the estimate \eqref{approx 2} holds.\cqfd

\end{document}